\documentclass[12pt,twoside]{amsart}
\usepackage{amssymb}

\nonstopmode

\numberwithin{equation}{section}
\newtheorem{theorem}{Theorem}[section]
\newtheorem{lemma}[theorem]{Lemma}
\newtheorem{proposition}[theorem]{Proposition}
\newtheorem{corollary}[theorem]{Corollary}

\newtheorem{claim}[theorem]{Claim}

\theoremstyle{definition}
\newtheorem{definition}[theorem]{Definition} 
\newtheorem{remark}[theorem]{Remark}
\newtheorem{example}[theorem]{Example}

\newcommand\lcm{\operatorname{lcm}}
\newcommand\length{\operatorname{length}}
\newcommand\pd{\operatorname{pd}}
\newcommand\codim{\operatorname{codim}}
\newcommand\coker{\operatorname{coker}}

\newcommand\Hom{\operatorname{Hom}}
\newcommand\Ext{\operatorname{Ext}}

\newcommand\rank{\operatorname{rank}}

\newcommand\depth{\operatorname{depth}}

\newcommand\im{\operatorname{im}}

\newcommand{\ER}{\operatorname{Ext}_R}
\newcommand{\ES}{\operatorname{Ext}_S}

\newcommand{\proj}[1]
{ \mathchoice
           { {\mathbb P}^{#1} }
           { {\mathbb P}^{#1} }
           { {\mathbb P}^{#1} }
           { {\mathbb P}^{#1} }
         }

\newcommand{\acm}{arithmetically Cohen-Macaulay }

\begin{document}
\title[Lifting Monomial Ideals]{Lifting Monomial Ideals}

\author[J.\ Migliore, U.\ Nagel]{J.\ Migliore$^*$, U.\ Nagel$^{**}$}

\address{Department of Mathematics,
Room 370 CCMB,
        University of Notre Dame, 
        Notre Dame, IN 46556, 
        USA}
\email{Juan.C.Migliore.1@nd.edu}

\address{Fachbereich Mathematik und Informatik, Universit\"at-Gesamthochschule
Paderborn, D--33095 Paderborn, Germany}
\email{uwen@uni-paderborn.de}

\date{July 6, 1999}
\thanks{$^*$ Partially supported  by the Department of Mathematics of the
  University of Paderborn\\
$^{**}$ Partially supported  by the Department of Mathematics of the
  University of Notre Dame }


\begin{abstract} We show how to lift any monomial ideal $J$ in $n$ variables
to a saturated ideal $I$ of the same codimension in $n+t$ variables.  We show
that $I$ has the same graded Betti numbers as $J$ and we show how to obtain
the matrices for the resolution of $I$.  The cohomology of $I$ is described. 
Making general choices for our lifting, we show that $I$ is the ideal of a
reduced union of linear varieties with singularities that are ``as small as
possible'' given the cohomological constraints.  The case where $J$ is
Artinian is the nicest.  In the case of curves we obtain stick figures for
$I$, and in the case of points we obtain certain $k$-configurations which we
can describe in a very precise way.
\end{abstract}


\maketitle
\tableofcontents


 \section{Introduction} \label{intro} 
Let $J \subset K[X_1,\dots,X_n]$ be the ideal of a subscheme, $W$, of
projective space $\proj{n-1}$.   An important question in general is to
determine what subschemes $V$ of $\proj{n}$ have $W$ as hyperplane section. 
More generally, we seek subschemes of $\proj{n+k}$ with $W$ as the
intersection with a general linear space of complementary dimension.  
Furthermore, one would like to find $V$ with nice properties, and to describe
the minimal free resolution of $V$.  If $W$ is an \acm subscheme of
$\proj{n-1}$ of dimension $\geq 1$, our construction should give the same
property for $V$.

It will be convenient to use new letters for the additional
variables.  Hence algebraically, in the hyperplane section case we seek a
saturated ideal $I \subset K[X_1,\dots,X_n,u_1]$ so that $J =
(I,u_1)/(u_1)$.  In this context, the problem makes sense even if $J$ is
not saturated (or even if $J$ is Artinian so there is no $W$).  In the more
general setting, we seek $I \subset K[X_1,\dots,X_n,u_1,\dots,u_t]$ so that $J
= (I ,u_1,\dots,u_t)/(u_1,\dots,u_t)$.  This is a special case of the
so-called ``lifting problem.''

An obvious solution is to simply view $J
\subset K[X_1,\dots,X_{n+1}]$ and consider the cone over $W$.  This has
exactly the same resolution as $J$, but a nasty singularity at the vertex
point.  One solution, in the case where $W$ is a curvilinear zeroscheme, can
be found in \cite{ballico-m} and \cite{walter1}, where they show that $V$ can
be taken to be even a smooth curve.  However, the curves obtained are almost
never arithmetically Cohen-Macaulay, and in \cite{walter2} it is shown that
``most of the time'' there does not even exist a smooth \acm curve with $W$ as
hyperplane section.

A better solution, from our point of view, is for the case of codimension two
varieties.  In this case the \acm property means that the minimal free
resolution of $J$ is a short exact sequence, so in this case the only map to
worry about is given by the Hilbert-Burch matrix (since the next map is
given by the maximal minors of the Hilbert-Burch matrix).  One can then add a
new variable, say $t$, to each entry of the matrix in a sufficiently general
way and obtain an \acm $V$ with $W$ as hyperplane section.  This was done in
\cite{CO} in the case where $W$ is a finite set of points in $\proj{2}$, and
they showed that if $W$ is of a certain (very general) form then $V$ will be a
smooth curve.

When the resolution is longer, however, there are more matrices in the
resolution and it is very difficult to add variables to each entry and still
preserve exactness (or even the property of being a complex).  Even finding
an ideal of the right codimension that restricts to $J$ (say by replacing each
occurrence of a variable in a generator of $J$ by a linear form involving the
new variables $u_1,\dots,u_t$) is difficult.  This was done in \cite{HTV} for
the case of codimension three arithmetically Gorenstein schemes, where the
original schemes were lifted to reduced irreducible ones, but as in the
codimension two case the whole resolution depends only on one matrix, so the
problem of preserving exactness came for free.

Failing smoothness or irreducibility, a desirable property for curves in
$\proj{3}$ is to be a so-called stick figure.  The study of such curves,
especially from the point of view of the Hilbert scheme, can be found in
\cite{hart-zeuthen}, where Hartshorne solved the long-standing Zeuthen
problem.  An extension of this property to codimension two varieties can be
found in \cite{BM5}.  In the case of dimension one, such curves were
constructed in \cite{RR} which were arithmetically Cohen-Macaulay, using a
method different from that in this paper.

The goal of this paper is to give a complete solution to the lifting problem 
in the case where $J$ is a monomial ideal and we add any number of
variables.  We lift to a reduced union of linear varieties (see below).  We
build on work of Hartshorne \cite{hart} and of Geramita, Gregory and Roberts
\cite{GGR}.  In those papers the authors began with an Artinian monomial
ideal $J$ and produced a ``lifted" ideal $I$ which they showed was the
saturated ideal of a reduced set of points.   Of course once the new ideal is
constructed, it is no longer monomial, so the process cannot be repeated. 
Our approach is to mimic the construction mentioned above, but to introduce
any number of variables at one time.  We make the procedure somewhat easier
to follow by introducing a matrix, called the ``lifting matrix," which
contains all the information of the lifting.  

In section 2 we describe our lifting method and we make a detailed  study of
the minimal free resolution of the ideal we obtain, as well as a study of the
cohomology.  The main tools are Taylor's free resolution of a monomial
ideal \cite{taylor} and the theorem of Buchsbaum and Eisenbud \cite{B-E} on
what makes a complex exact. A lemma of Buchsbaum and Eisenbud on lifting
\cite{BE-lifting} is also useful.  The main result of the section is
Proposition \ref{same-resolution}, where we show that the graded Betti
numbers of the lifted ideal are exactly the same as those of $J$, and we
describe how the maps in the resolution are lifted from those of the minimal
free resolution of $J$.  As a consequence we get that $I$ is a saturated
ideal, and we obtain a great deal of cohomological information.  For
instance, the depth of the lifted ideal is large, and the
cohomology of $I$ is a ``lifting'' (in a precise sense) of that of $J$. This
section requires almost nothing about the lifting matrix, and as a result
gives almost nothing (except dimension) in the way of ``nice'' geometric
properties of the scheme defined by the lifted matrix.  We even observe in
Remark \ref{more-general} that the lifting matrix can be chosen in such a way
that the new ideal is not, strictly speaking, a lifting but does still have
the cohomological properties just described.  This fact is utilized in 
section 4.   We call such an ideal a ``pseudo-lifting.''

In section 3 we assume that the entries of the lifting matrix are
sufficiently general.  We define an extension of the notion of stick figures
to any dimension or codimension.  The strongest condition is what we call a
``generalized stick figure.''  We prove in Theorem \ref{lift-to-stick} that
our schemes $V$ obtained by using this general lifting matrix are not only
reduced, but in fact the union of the components of any given dimension are
generalized stick figures away from $W$.   As a corollary, if $J$ is Artinian
then $V$ is an arithmetically Cohen-Macaulay generalized stick figure. 
Furthermore, we can construct in this way an \acm generalized stick figure
corresponding to any allowable Hilbert function (Corollary
\ref{differentiable}).  Moreover, any monomial ideal is the initial ideal of
a radical ideal, and any Artinian monomial ideal is the initial ideal of a
radical ideal defining an \acm generalized stick figure.

In section 4 we give a more detailed description of the configurations of
linear varieties obtained by lifting Artinian monomial ideals, both in the
strict sense and in the more general sense of pseudo-liftings mentioned
above.  It turns out that the case of lex-segment Artinian monomial ideals
gives a particularly nice special case.  They produce certain so-called
``$k$-configurations,'' while the other Artinian monomial ideals produce only
so-called ``weak
$k$-configurations.''  Using our approach we suggest an extension of
$k$-configurations to higher dimension.

It is a pleasure to acknowledge the many contributions of Robin Hartshorne
in the area of this paper and in the broader field of Algebraic Geometry.
We dedicate this paper to him on the occasion of his sixtieth birthday.


\section{Lifting Monomial Ideals: Most General Case}

Let $k$ be an infinite field and let $S = K[X_1,\dots,X_n]$ and
$R=K[X_1,\dots,X_n,u_1,\dots,u_t]$.  Our
techniques also often work over a finite field: see Remark
\ref{matrix-condition}.   

We shall be interested in the general idea of ``lifting'' a monomial ideal. 
More generally, we have

\begin{definition}\label{lift-def}
Let $R$ be a ring and let $u_1,\dots,u_t$ be elements of $R$ such that 
$\{u_1,\dots,u_t\}$ forms an $R$-regular sequence.  Let $S =
R/(u_1,\dots,u_t)$.  Let $B$ be an $S$-module and let $A$ be an $R$-module. 
Then we say $A$ is a {\em $t$-lifting} of $B$ to $R$ if $\{u_1,\dots,u_t\}$
is an $A$-regular sequence and $A/(u_1,\dots,u_t)A \cong B$.  When $t=1$ we
will sometimes just refer to $A$ as a {\em lifting} of $B$.
\end{definition}

The following lemma of Buchsbaum and
Eisenbud \cite{BE-lifting} will be very useful.

\begin{lemma}\label{BE-lifting-lemma}
Let $R$ be a ring, $x \in R$, and $S = R/(x)$.  Let $B$ be an $S$-module, and
let 
\[
{\mathcal F}: F_2 \stackrel{\phi_2}{\longrightarrow} F_1
\stackrel{\phi_1}{\longrightarrow} F_0
\]
be an exact sequence of $S$-modules with $\coker \phi_1 \cong B$.  Suppose
that
\[
\Gamma:G_2 \stackrel{\psi_2}{\longrightarrow} G_1
\stackrel{\psi_1}{\longrightarrow} G_0
\]
is a complex of $R$-modules such that 
\begin{itemize}
\item[(1)] The element $x$ is a non-zero divisor on each $G_i$,
\item[(2)] $G_i \otimes_R S = F_i$, and
\item[(3)] $\psi_i \otimes_R S = \phi_i$.
\end{itemize}
Then $A = \coker \psi_1$ is a lifting of $B$ to $R$.
\end{lemma}

In Definition \ref{lift-def}, notice that if $A$ is an ideal in a
polynomial ring then $A$ is a lifting of $B$ if and only if $x$ is not a 
zero-divisor on $R/A$ and $(A,x)/(x) \cong B$.  To set notation we rewrite
Definition \ref{lift-def} for the case of ideals.

\begin{definition}\label{t-lifting-def}
Let $R= K[X_1,\dots,X_n,u_1,\dots,u_t]$ and $S = K[X_1,\dots,X_n]$.  Let $I
\subset R$ and $J \subset S$ be homogeneous ideals.  Then we say $I$ is a {\em
$t$-lifting} of $J$ to $R$ (or when $R$ is understood, simply a {\em
$t$-lifting} of $J$) if $(u_1,\dots,u_t)$ is a regular sequence
on $R/I$ and $(I,u_1,\dots,u_t)/(u_1,\dots,u_t) \cong J$.  We say that $I$ is
a  {\em reduced $t$-lifting} of $J$ if it is a $t$-lifting and if furthermore
$I$ is a radical ideal in $R$.
\end{definition}

In this section we show how to lift a monomial ideal, adding any number of
variables.  Later we will make ``general'' selections so that the lifted
ideals will be the saturated ideals of projective subschemes which are not
only reduced but in fact {\em generalized stick figures} (see Definition
\ref{gen-stick-fig}).

In \cite{GGR} the authors show, following an idea of Hartshorne \cite{hart},
how to get a reduced 1-lifting of a monomial ideal.  The ideal produced in
this way is no longer a monomial ideal, so their construction cannot be
repeated to produce higher lifting.  Furthermore, although they do not
explicitly say so, their monomial ideals seem to be Artinian, since they
describe their lifted ideal geometrically as the ideal of polynomials
vanishing on a certain set
$\bar M$, and they note that $u_1$ (in our notation) is not a zero divisor
since it does not vanish at any element of $\bar M$, suggesting that $\bar M$
is a finite set.  In this section we will show how to $t$-lift a monomial
ideal by in fact lifting the whole resolution.  This gives a great deal of
information about the lifted ideal.

Our approach to the lifting follows that of \cite{GGR}, and we now recall some
notation introduced there.  Let ${\mathbb N} = \{ 0,1,2,\dots \}$ and if
$\alpha = (a_1,\dots,a_n) \in {\mathbb N}^n$, we let
$X^\alpha = X_1^{a_1}\cdots X_n^{a_n}$.  Letting $P$ denote the set of
monomials in $S$ (including $1$), then this gives a bijection between $P$ and
${\mathbb N}^n$.  The set ${\mathbb N}^n$ may be partially ordered by
$(a_1,a_2,\dots,a_n) \leq (b_1,b_2,\dots,b_n)$ if and only if $a_i \leq b_i$
for all $i$.  This partially ordering translates, via the bijection above, to
divisibility of monomials in $P$.

For each variable $X_j$, $1 \leq j \leq n$, choose infinitely many 
linear forms $L_{j,i} \in$ \break $ K[X_j,u_1,\dots,u_t]$ ($i =1,2,\dots$). 
In subsequent sections we will impose other conditions so that each $L_{j,i}$
is chosen sufficiently generically with respect to $L_{j,1},\dots,L_{j,i-1}$,
but for this section we do not even assume that no $L_{j,i}$ is a scalar
multiple of an $L_{j,k}$.  We only assume that the coefficient of $X_j$ in
$L_{j,i}$ is not zero.  For any given example, only finitely many
$L_{j,i}$ need be chosen.  In this case the matrix $A = [L_{j,i}]$ will be
called the {\em lifting matrix}:
\[
A = 
\left [
\begin{array}{ccccccccc}
L_{1,1} & L_{1,2} & L_{1,3} & \dots \\
L_{2,1} & L_{2,2} & L_{2,3} & \dots \\
\vdots & \vdots & \vdots \\
L_{n,1} & L_{n,2} & L_{n,3} & \dots 
\end{array}
\right ].
\]

To each monomial $m \in P$, if $m = \prod_{j=1}^n X_j^{a_j}$, we associate the
homogeneous polynomial 
\begin{equation}\label{lift-m-def}
\bar m = \prod_{j=1}^n \left ( \prod_{i=1}^{a_j} L_{j,i} \right ) \in R.
\end{equation}
If $J = (m_1,\dots,m_r)$ is a monomial ideal in $S$ then we denote by $I$ the
ideal $(\bar m_1,\dots,\bar m_r ) \subset R$.  The crucial aspect of this
construction is that at each occurrence of some power $a$ of $X_j$ in any
monomial, we always take the product $L_{j,1}\cdots L_{j,a}$ of the {\em first
$a$ entries} of the $j$-th row (corresponding to $X_j$) of the lifting
matrix.  This redundancy is what makes the construction work.

\begin{remark} \label{geometrically}
Geometrically the choice of the $L_{j,i}$ above amounts to viewing
$\proj{n-1} \subset \proj{n+t-1}$ and choosing  hyperplanes in
$\proj{n+t-1}$ containing the hyperplane in $\proj{n-1}$ defined by $X_j$. 
\end{remark}

\begin{remark} \label{polarization}
The construction above is a generalization, as mentioned, of the
construction of \cite{GGR}, which gave the case $t=1$.  We would like to also
point out that the unpublished thesis of Schwartau \cite{schwartau} contains a
construction which is also similar.  That is, his notion of {\em
polarization} of a monomial ideal is obtained by replacing the repeated
variables in a monomial ideal by {\em new variables} instead of by different
linear forms in one or $t$ new variables.  His conclusions are
different from ours, although some of his preparatory lemmas are useful and
are quoted below.
\end{remark}

We now show how to lift the minimal free resolution of a monomial ideal. The
first step is to use Taylor's (not necessarily minimal) free resolution
\cite{taylor}.  Recall that the entries of the matrices of such a resolution
are themselves monomials (with a suitable choice of basis).  These entries
will be replaced by successive copies of the $L_{j,i}$ in ``almost'' the
analogous way to what was done for the generators.

\begin{proposition} \label{same-resolution}
Let $J = (m_1,\dots,m_r) \subset S$ be a monomial ideal and let $I = (\bar
m_1,\dots,\bar m_r)$ be the ideal described above.  Consider a (minimal) free
$S$-resolution
\[
0 \rightarrow F_p \stackrel{\phi_p}{\longrightarrow} F_{p-1}
\stackrel{\phi_{p-1}}{\longrightarrow} \cdots
\stackrel{\phi_2}{\longrightarrow} F_1 \stackrel{\phi_1}{\longrightarrow} J
\rightarrow 0
\]
of $J$.  Then $I$ has a (minimal) free $R$-resolution
\begin{equation}\label{lifted-resolution}
0 \rightarrow \bar F_p \stackrel{\bar \phi_p}{\longrightarrow} \bar F_{p-1}
\stackrel{\bar \phi_{p-1}}{\longrightarrow} \cdots
\stackrel{\bar \phi_2}{\longrightarrow} \bar F_1
\stackrel{\bar \phi_1}{\longrightarrow} I \rightarrow 0
\end{equation}
where $\bar F_i$ is a ``lifting'' of $F_i$ in the obvious way and the maps
$\bar \phi_i$ are ``liftings,'' explicitly obtained from the $\phi_i$ as
described below in the proof.
\end{proposition}

\begin{proof}
For the convenience of the reader, we first briefly recall from
\cite{eisenbud} the relevant facts about Taylor's resolution.  Let
$m_1,\dots,m_r$ be monomials in $S$.  Taylor's free resolution of the ideal
$J = (m_1,\dots,m_r)$ has the form
\[
0 \rightarrow F_r \stackrel{d_r}{\longrightarrow} \cdots
\stackrel{d_2}{\longrightarrow} F_1 
\stackrel{d_1}{\longrightarrow} J \rightarrow 0
\] 
where the free modules and maps are defined as follows.  Let $F_s$ be the free
module on basis elements $e_A$, where $A$ is a subset of length $s$ of $\{
1,\dots,r \}$.  Set 
\[
m_A = \lcm \{ m_i | i \in A \}.
\]
For each pair $A,B$ such that $A$ has $s$ elements and $B$ has $s-1$
elements, let $A = \{ i_1,\dots,i_s \}$ and suppose that $i_1 < \cdots <
i_s$.  Define 
\[
c_{A,B} = 
\left \{
\begin{array}{ll}
0 & \hbox{if } B \not\subset A \\
(-1)^k m_A/m_B & \hbox{if } A = B \cup \{i_k \} \hbox{ for some } k.
\end{array}
\right.
\]
Then we define $d_s : F_s \rightarrow F_{s-1}$ by sending $e_A$ to $\sum_B
c_{A,B} e_B$.
Clearly Taylor's resolution does not, in general, have the same length as the
minimal free resolution of $J$ since the length is equal to the number of
generators of $J$ and, for example, $F_r$ has rank 1.

Now, a key observation is given in \cite{eisenbud} at the end of
Exercise 17.11 (page 439): the exact same construction can be used in a much
more general setting to produce at least a complex.  In particular, in our
situation, replacing the monomial ideal $J = (m_1,\dots,m_r) \subset S$ by the
ideal $I = (\bar m_1,\dots,\bar m_r) \subset R$ , we see that we 
at least have a complex.
Furthermore, since $m_A$ is also a monomial it can be lifted as above, and one
can immediately check that we have 
\begin{equation}\label{lcm-fact}
\bar m_A = \lcm \{ \bar m_i | i \in A \}.
\end{equation}
Note that we do not necessarily have the same kind of lifting for the
$c_{A,B}$; it is a lifting, but the products of  $L_{j,i}$ do not necessarily
begin with $L_{j,1}$, as we indicated above (see Example
\ref{clarifying-example}).  However, this does not matter.
We conclude that we have a complex
\[
0 \rightarrow \bar F_r \stackrel{\bar d_r}{\longrightarrow} \cdots
\stackrel{\bar d_2}{\longrightarrow} \bar F_1 \stackrel{\bar
d_1}{\longrightarrow} I
\rightarrow 0
\]
where the $\bar d_k$ restrict to the $d_k$.

It remains to check that this complex is in fact a resolution for $I$.  By
the Buchsbaum-Eisenbud exactness criterion, we have to show that (i) $\rank
\bar d_{k+1} + \rank \bar d_k = \rank \bar F_k$ for all $k$ and (ii) the
ideal of $(\rank \bar d_k )$-minors of $\bar d_k$ contains a regular sequence
of length $k$, or is equal to $R$.  Both of these follow immediately from the
fact that the restriction of the $\bar d_k$ is $d_k$ for each $k$, and we
know that the restriction is Taylor's resolution and hence satisfies these
two properties.

Finally, an entry of one of the $d_k$ is 1 if and only if the
corresponding lifted matrix has a 1 in the same position.  Hence also the
minimal free resolutions agree, as claimed.
\end{proof}

\begin{example}\label{clarifying-example}
Let $n=3$ and $t=2$.  Consider the ideal $J = (X_1^2X_2,X_2^2X_3,X_3^2X_1)$. 
This has the minimal free $S$-resolution
\[
0 \rightarrow S(-6) \stackrel{\phi_3}{\longrightarrow} S(-5)^3
\stackrel{\phi_2}{\longrightarrow} S(-3)^3 \stackrel{\phi_1}{\longrightarrow}
J \rightarrow 0
\]
where
\[
\begin{array}{ccc}

\phi_1 = \left [ 
\begin{array}{ccc}
X_1^2X_2 & X_2^2X_3 & X_3^2X_1
\end{array}
 \right ] \\ \\

\phi_2 = \left [
\begin{array}{ccc}
-X_3^2 & 0 & -X_2X_3 \\
0 & -X_1X_3 & X_1^2 \\
X_1X_2 & X_2^2 & 0
\end{array}
\right ]
\\ \\

\phi_3 = \left [
\begin{array}{c}
-X_2 \\
X_1 \\
X_3
\end{array}
\right ]
\end{array}
\]
Notice that in this case Taylor's resolution is minimal.  In order to lift, we
choose linear forms for each variable as follows:
\[
\begin{array}{c}
X_1 \ : \ L_{1,1}(X_1, u_1, u_2), L_{1,2}(X_1, u_1, u_2), \dots \\
X_2 \ : \ L_{2,1}(X_2, u_1, u_2), L_{2,2}(X_2, u_1, u_2), \dots \\
X_3 \ : \ L_{3,1}(X_3, u_1, u_2), L_{3,2}(X_3, u_1, u_2), \dots 
\end{array}
\]
Then we set $I = (L_{1,1}L_{1,2}L_{2,1}, L_{2,1}L_{2,2}L_{3,1},
L_{3,1}L_{3,2}L_{1,1})$ and its minimal free resolution has the form
\[
0 \rightarrow R(-6) \stackrel{\bar \phi_3}{\longrightarrow} R(-5)^3
\stackrel{\bar \phi_2}{\longrightarrow} R(-3)^3
\stackrel{\bar \phi_1}{\longrightarrow} I \rightarrow 0
\]
where
\[
\begin{array}{ccc}

\bar \phi_1 = \left [ 
\begin{array}{ccc}
L_{1,1}L_{1,2}L_{2,1} & L_{2,1}L_{2,2}L_{3,1} & L_{3,1}L_{3,2}L_{1,1}
\end{array}
 \right ] \\ \\

\bar \phi_2 = \left [
\begin{array}{ccc}
-L_{3,1}L_{3,2} & 0 & -L_{2,2}L_{3,1} \\
0 & -L_{1,1}L_{3,2} & L_{1,1}L_{1,2} \\
L_{1,2}L_{2,1} & L_{2,1}L_{2,2} & 0
\end{array}
\right ]
\\ \\

\bar \phi_3 = \left [
\begin{array}{c}
-L_{2,2} \\
L_{1,2} \\
L_{3,2}
\end{array}
\right ]
\end{array}
\]
Notice, for example, that the entries of $\bar \phi_3$ are all of the form
$L_{j,2}$ rather than $L_{j,1}$.
\end{example}

\begin{corollary}\label{depth-cor}
$\depth R/I = \depth S/J + t$.  
In particular, $\depth R/I \geq t$.
\end{corollary}

\begin{proof}
\[
\begin{array}{rcl}
\depth R/I & = & n+t - \pd_R R/I \\
& = & n+t - \pd_S S/J \\
& = & n+t - [n - \depth S/J]
\end{array}
\]
\end{proof}

For the rest of this paper we will be interested in the projective
subschemes defined by $I$ and $J$:

\begin{quotation}
{\em We denote by $V$ the scheme in
$\proj{n+t-1}$ defined by $I$, and by $W$ the scheme in $\proj{n-1}$ defined
by $J$.  If this latter is empty, i.e.\ if $S/J$ is Artinian, then for the
purposes below we formally define $\codim W = n$ and $\deg W = \length S/J$.}
\end{quotation}

\begin{lemma} \label{codim}
$V$ has the same codimension in $\proj{n+t-1}$ that $W$ has in $\proj{n-1}$.
\end{lemma}

\begin{proof}
Suppose that the codimension of $W$ is $c$.  Clearly $W$ is the intersection
in $\proj{n+t-1}$ of $V$ with the codimension $t$ linear space defined by
$u_1 = \dots = u_t = 0$.  Hence $\codim V \geq c$.  So we only have to prove
$\codim V \leq c$.

Since $J$ is a monomial ideal, all associated primes are of the form
$(X_{i_1},\dots,X_{i_k})$.  By hypothesis, then, there exist
$X_{i_1},\dots,X_{i_c}$ such that every element of $J$ is in the ideal
$(X_{i_1},\dots,X_{i_c})$.  By the construction of $I$ it is clear that every
element of $I$ is in the ideal $(L_{i_1,1},\dots,L_{i_c,1})$.  Hence $\codim
V \leq c$ as claimed.
\end{proof}

\begin{corollary} \label{some-facts}
\begin{itemize}
\item[(i)] The ideal $I$ is saturated.
\item[(ii)] $S/J$ is Cohen-Macaulay (including the case where it is
Artinian) if and only if $R/I$ is Cohen-Macaulay.
\item[(iii)] $(I,u_1,\dots,u_t)/(u_1,\dots,u_t) \cong J$.
\item[(iv)]  $\deg V = \deg W$.
\item[(v)] $(u_1,\dots,u_t)$ is a regular sequence on $A = R/I$.
\end{itemize}
Combining (iii) and (v), we get that $I$ is a $t$-lifting of $J$.
\end{corollary}

\begin{proof}
(i) is immediate from Corollary \ref{depth-cor}.  For (ii), $\Rightarrow$
follows from Corollary \ref{depth-cor} and Lemma \ref{codim} (the converse is
immediate).  (iii) is obvious.  (iv) follows from Proposition
\ref{same-resolution} and a computation of the Hilbert polynomial. 

To prove (v) we use induction on $t$.  The case $t=1$ follows from
Proposition \ref{same-resolution} and Lemma \ref{BE-lifting-lemma}.  Let
$I_i$ ($1 \leq i \leq t$) denote the ideal in $R_i :=
K[X_1,\dots,X_n,u_1,\dots,u_i]$ which lifts $J$ to that ring via our
construction.  Note $I = I_t$; in this case we continue to refer to the ideal
as $I$. By induction we may assume that $(u_1,\dots,u_{t-1})$ is a regular
sequence on $R_{t-1}/I_{t-1}$.  Then by Proposition \ref{same-resolution} we
see that setting $R = K[X_1,\dots,X_n,u_1,\dots,u_t]$, $x = u_t$, $S =
K[X_1,\dots,X_n,u_1,\dots,u_{t-1}]$ and $B = S/I_{t-1}$, all the hypotheses of
Lemma \ref{BE-lifting-lemma} are satisfied, and hence $R/I$ is a lifting of
$R_{t-1}/I_{t-1}$ to $R$.  The result follows immediately.
\end{proof}

The following result from \cite{bruns-herzog} (Lemma 3.1.16, p.\ 94) will be
useful.  Here we give the graded version.

\begin{lemma}\label{bruns-herzog-lemma}
Let $R$ be a graded ring, and let $M$ and $N$ be graded $R$-modules.  If $x$
is a homogeneous $R$- and $M$-regular element with $x \cdot N = 0$, then 
\[
\Ext_R^{i+1} (N,M)(-\deg x) \cong \Ext^i_{R/(x)} (N,M/xM)
\]
for all $i \geq 0$.
\end{lemma}

\begin{proposition} \label{t-lift-ext}
If $\Ext^i_S (S/J,S) \neq 0$ for some $i \in {\mathbb Z}$ then $\Ext_R^i
(R/I,R)$ is a $t$-lifting of $\Ext_S^i (S/J,S)$.  In particular, $\depth
\Ext_R^i (R/I,R) = t + \depth^i_S (S/J,S)$.
\end{proposition}

\begin{proof}
We adopt the following notation:
\[
\begin{array}{rcl}
( - )^* & = & \Hom_R (-,R) \hbox{ (dualizing with respect to $R$)} \\
( - )^\vee & = & \Hom_S (-,S) \hbox{ (dualizing with respect to $S$)}
\end{array}
\]
Let $t=1$ and put $u = u_1$.  Let
\[
{\mathbb F}_{\bullet} : \ \ \ \cdots \rightarrow F_{i+1}
\stackrel{\varphi_{i+1}}{\longrightarrow} F_i
\stackrel{\varphi_i}{\longrightarrow} \cdots
\stackrel{\varphi_1}{\longrightarrow} S \rightarrow S/J \rightarrow 0
\]
be the Taylor resolution of $S/J$.  From the proof of Proposition
\ref{same-resolution} we know that $R/I$ has a resolution
\[
{\bar \mathbb F}_{\bullet} : \ \ \ \cdots \rightarrow \bar F_{i+1}
\stackrel{\bar \varphi_{i+1}}{\longrightarrow} \bar F_i
\stackrel{\bar \varphi_i}{\longrightarrow} \cdots
\stackrel{\bar \varphi_1}{\longrightarrow} R \rightarrow R/I \rightarrow 0
\]
such that $\bar {\mathbb F}_\bullet \otimes_R S = {\mathbb F}_\bullet$.  In
other words, we have a short exact sequence of complexes
\[
0 \rightarrow \bar {\mathbb F}_\bullet (-1) \stackrel{u}{\longrightarrow}
\bar {\mathbb F}_\bullet \rightarrow {\mathbb F}_\bullet \rightarrow 0.
\]
Dualizing we obtain the short exact sequences of complexes
\[
0 \rightarrow \bar {\mathbb F}_\bullet^* (-1) \stackrel{u}{\longrightarrow} 
\bar {\mathbb F}_\bullet^* \rightarrow \Ext_R^1 ({\mathbb F}_\bullet ,R)(-1)
\rightarrow 0,
\]
since
\[
0 \rightarrow \bar {F}_i (-1) \stackrel{u}{\longrightarrow} \bar
{F}_i \rightarrow {F}_i \rightarrow 0
\]
induces
\[
\begin{array}{ccccccccccccccc}
0 & \! \! \! \rightarrow & \! \! \!  \Hom_R (F_i ,R) & \! \! \!  \rightarrow
& \! \! \!  \bar F_i^* & \! \! \! 
\rightarrow & \! \! \!  \bar F_i^* (1) & \! \! \!  \rightarrow & \! \! \! 
\Ext_R^1(F_i,R) & \! \! \! 
\rightarrow & \! \! \! 
\Ext_R^1 (\bar F_i,R).\\
&&  \! \! \! || &&&&&&&&  \! \! \!  || \\
&&  \! \! \! 0  &&&&&&&&  \! \! \!  0
\end{array}
\]
Lemma \ref{bruns-herzog-lemma} gives 
\[
\Ext_R^1 ({\mathbb F}_\bullet ,R)(-1) \cong \Hom_S ({\mathbb F}_\bullet ,S) =
{\mathbb F}_\bullet^\vee.
\]
Thus we get the short exact sequence of complexes
\begin{equation} \label{seq-of-complexes}
0 \rightarrow \bar {\mathbb F}_\bullet^* (-1) \stackrel{u}{\longrightarrow}
\bar {\mathbb F}_\bullet^* \rightarrow {\mathbb F}_\bullet^\vee \rightarrow 0.
\end{equation}
From the exact sequence
\[
0 \rightarrow R/I(-1) \stackrel{u}{\longrightarrow} R/I \rightarrow S/J
\rightarrow 0
\]
and using Lemma \ref{bruns-herzog-lemma} and twisting by $-1$ we get the
induced long exact homology sequence
\[
\begin{array}{c}
\cdots \rightarrow \Ext_S^{i-1} (S/J,S) \rightarrow \Ext^i_R (R/I,R)(-1)
\stackrel{u}{\longrightarrow} \Ext^i_R(R/I,R) \rightarrow \\ \\ 
 \Ext^i_S
(S/J,S) \rightarrow \Ext_R^{i+1} (R/I,R)(-1) \stackrel{u}{\longrightarrow}
\cdots
\end{array}
\]
\begin{claim} 
The multiplication map $\Ext_R^i (R/I,R)(-1) \stackrel{u}{\longrightarrow}
\Ext_R^i (R/I,R)$ is injective for all $i \geq 1$.
\end{claim}

The claim immediately proves the assertion of this proposition since it
allows us to break the homology sequence into short exact sequences of the
form we want.

In order to show the claim we consider the commutative diagram
\[
\begin{array}{cccccccccc}
&&&& \! \!  0 \\
&&&& \! \!  \downarrow \\ 
0 & \! \! \rightarrow & \! \!  \im (\bar \varphi_i^* )(-1) & \! \! 
\rightarrow & \! \!  \bar F_i^* (-1) & \! \!  \rightarrow & \! \!  \bar
F_i^*/\im(\bar \varphi_i^*)(-1) & \! \!  \rightarrow & \! \!  0 \\  
&& \! \!  \phantom{u} \downarrow u && \! \!  \phantom{u} \downarrow u && \!
\!  \phantom{u} \downarrow u \\ 
0 & \! \!  \rightarrow & \! \!  \im (\bar \varphi_i^* ) & \! \!  \rightarrow
& \! \!  \bar F_i^*  & \! \!  \rightarrow & \! \!  \bar F_i^* /\im(\bar
\varphi_i^*) & \! \!  \rightarrow & \! \!  0 \\  
&&&& \! \!  \downarrow \\
&&&& \! \!  F_i^\vee \\
&&&& \! \!  \downarrow \\
&&&& \! \!  0
\end{array}
\]
where the center column comes from (\ref{seq-of-complexes}).  Denote by
$\alpha$ the multiplication on the left-hand side.  Then the snake lemma
provides a homomorphism $\tau :
\coker
\alpha
\rightarrow F_i^\vee$.  Since $\bar F_\bullet \otimes_R S = F_\bullet$ we
have $\coker
\alpha \cong
\im (\varphi_i^\vee)$.  Thus $\tau$ is induced from the embedding
$\im(\varphi_i^\vee) \hookrightarrow F_i^\vee$, i.e.\ $\tau$ is injective. 
Therefore the snake lemma shows that
\[
\bar F_i^* /\im(\bar \varphi_i^*)(-1) \stackrel{u}{\longrightarrow} \bar
F_i^*/\im(\bar \varphi_i^*)
\]
is injective.  Since $\Ext_R^i (R/I,R)$ is a submodule of $\bar F_i^*/\im
(\bar
\varphi_i^*)$, our claim is proved.

Now let $t>1$.  In the argument above, we have not used the specific lifting
from $J$ to $I$, but only the fact that a resolution of $J$ lifts to a
resolution of $I$.  Hence our assertion follows by induction on $t$.
\end{proof}

\begin{corollary}\label{acm-or-not}
$V$ is either arithmetically Cohen-Macaulay (if $S/J$ is Cohen-Macaulay) or
else it fails to be both locally Cohen-Macaulay and equidimensional.
\end{corollary}

\begin{proof}
We saw in Corollary \ref{some-facts} that $R/I$ is Cohen-Macaulay if and only
if $S/J$ is Cohen-Macaulay.  Assume, then, that $S/J$ is not Cohen-Macaulay. 
Suppose $V$ has codimension $c$; note that $c \leq n-1$.  It is enough to show
that $\ER^{i} (R/I,R)$ does not have finite length for some $i$ in the range
$c+1 \leq i \leq n$.  We have by assumption that $\ES^i (S/J,S)$ is non-zero
for some $i$ in the range $c+1 \leq i \leq n$.  Hence the result follows
immediately from Proposition \ref{t-lift-ext}.
\end{proof}

\begin{example}
Let $n=3$, $t=1$ and consider the ideal $J$ in $K[X_1,X_2,X_3]$ defined by
$(X_1^2, X_1X_2, X_1X_3, X_2^2,X_2X_3)$.  Notice that $J$ is not saturated,
but that its saturation $\tilde J$ is just $(X_1,X_2)$, hence $S/\tilde J$ is
Cohen-Macaulay but $J$ has an irrelevant primary component.  One can check
that the lifted ideal $I \subset K[X_1,X_2,X_3,u_1]$ is the saturated ideal
of the union of a line $\lambda \subset \proj{3}$ (defined by $X_1=X_2=0$)
and two points.  This is locally Cohen-Macaulay but not equidimensional.  Its
top dimensional part is arithmetically Cohen-Macaulay.

Now let $n=4$ and $t=1$ and consider the ideal $J$ in $K[X_1,X_2,X_3,X_4]$
defined by $(X_1X_3,X_1X_4,X_2X_3,X_2X_4)$.  $J$ is the saturated ideal of
the disjoint union of two lines in $\proj{3}$, hence $S/J$ is not
Cohen-Macaulay.  One can check that the lifted ideal $I \subset
K[X_1,X_2,X_3,X_4,u_1]$ is the saturated ideal of the union of two planes in
$\proj{4}$ meeting at a single point, which is equidimensional but not
locally Cohen-Macaulay.

These examples also illustrate Lemma \ref{decompose-I} and Theorem
\ref{lift-to-stick} below.
\end{example}

We conclude this section with some preparatory facts about
monomial ideals and liftings.  We will use them in the following section,
where we make some generality assumptions on our lifting matrix, but they
hold in the generality of the current section.

\begin{lemma}\label{schwartau-lemma}
Let 
\[
\begin{array}{rcl}
I^{(1)} & = & (m_1^{(1)},\dots,m_{n_1}^{(1)} ) \\
& \vdots \\
I^{(\ell)} & = & (m_1^{(\ell)},\dots,m_{n_\ell}^{(\ell)} )
\end{array}
\]
be monomial ideals in $S$.  Then $I^{(1)} \cap \dots \cap I^{(\ell)}$ is the
ideal generated by 
\[
\{ \lcm (m^{(1)}_{i_1} ,\dots, m_{i_\ell}^{(\ell)} ) \}
\]
where the indices lie in the ranges
\[
\begin{array}{c}
1 \leq i_1 \leq n_1 \\
\vdots \\
1 \leq i_\ell \leq n_\ell
\end{array}
\]
\end{lemma}

\begin{proof}
This is \cite{schwartau} Lemma 85.
\end{proof}

\begin{lemma}
Let $J_1, J_2 \subset S$ be monomial ideals and fix a lifting matrix $A$.
For any monomial ideal $J$, denote by $\bar J$ the lifting of $J$ by $A$. 
Then
\begin{itemize}
\item[(i)] $J_1 \subset J_2$ if and only if $ \bar I_1 \subset \bar I_2$.
\item[(ii)] $\overline {J_1 \cap J_2} = \bar J_1 \cap \bar J_2$.
\end{itemize}
\end{lemma}

\begin{proof}
Part (i) is clear.  For part (ii), the inclusion $\subseteq$ follows from
part (i).  Equality will come by showing
that the Hilbert functions are the same.  Consider the sequences
\[
\begin{array}{c}
0 \rightarrow J_1 \cap J_2 \rightarrow J_1 \oplus J_2 \rightarrow J_1 + J_2
\rightarrow 0 \\
0 \rightarrow \bar J_1 \cap \bar J_2 \rightarrow \bar J_1 \oplus \bar J_2
\rightarrow \bar J_1 + \bar J_2 \rightarrow 0
\end{array}
\]
Note that $\overline{J_1 + J_2} = \bar J_1 + \bar J_2$.  Thanks to Corollary
\ref{depth-cor} we then have the following calculation.  (We use the notation
$\Delta^t h_{R/I} (x)$ for the $t$-th difference of the Hilbert function of
$R/I$.  This is standard notation, but in any case see the discussion at the
end of section 3.)
\[
\begin{array}{rcl}
\Delta^t h_{R/ \overline{J_1 \cap J_2}} (x) & = & h_{S/J_1 \cap J_2} (x) \\
& = & h_{S/J_1} (x) + h_{S/J_2} (x) - h_{S/J_1 + J_2} (x) \\
& = & \Delta^t h_{R/\bar J_1} (x) + \Delta^t h_{R/\bar J_2} (x)
- \Delta^t h_{R/\overline{J_1 +J_2}} (x) \\
& = & \Delta^t h_{R/\bar J_1} (x) + \Delta^t h_{R/\bar J_2} (x) 
- \Delta^t h_{R/\bar J_1 + \bar J_2} (x) \\
& = & \Delta^t h_{R/\bar J_1 \cap \bar J_2 }(x)
\end{array}
\]
Hence the Hilbert functions agree.
\end{proof}

\begin{corollary} \label{decompose-I}
Let $J \subset S$ be a monomial ideal and let $I$ be the lifting of $J$
using some lifting matrix $A$.  Suppose that the primary decomposition of $J$
is
\[
J = Q_1 \cap \dots \cap Q_r.
\]
Then 
\[
I = \bar Q_1 \cap \dots \cap \bar Q_r
\]
where $\bar Q_i$ is the ideal generated by the liftings of the generators of
$Q_i$.
\end{corollary}

\begin{remark} \label{primary-decomp-of-mon-id}
It is known (cf.\ \cite{eisenbud} Exercise 3.8 or \cite{schwartau}, proof of
Theorem 91) that if $J \subset S$ is a monomial ideal then we can write $J =
Q_1 \cap \dots \cap Q_r$ where each $Q_i$ is a complete intersection of the
form $(X_{i_1}^{a_1} ,\dots, X_{i_p}^{a_p} )$ with $a_j \geq 1$ for all $j$,
$1 \leq j \leq p$.
\end{remark}

\begin{corollary}\label{unmixed}
With the convention that the empty set is equidimensional of dimension $-1$,
we have
\begin{itemize}
\item[(i)] $V$ is equidimensional if and only if $W$ is equidimensional.
\item[(ii)] If $W$ is equidimensional then $V$ is either \acm or not even
locally Cohen-Macaulay.
\end{itemize}
\end{corollary}

\begin{proof}
Part (i) follows from Corollary \ref{decompose-I} and Remark
\ref{primary-decomp-of-mon-id}.  Part (ii) also uses Corollary
\ref{acm-or-not} and its proof.
\end{proof}

\begin{remark}\label{more-general}
In defining our lifting matrix $A$, we said that in this section we assume
almost nothing about the linear forms which are its entries.  We required
only that the linear forms $L_{j,i}$ from the $j$-th row were elements of the
ring $K[X_j,u_1,\dots,u_t]$.  But there is another direction that we can go
with the techniques of this section.  Let $J$ be a monomial ideal in $S =
K[X_1,\dots,X_n]$ and let $A$ be an $n \times r$ matrix of linear forms
so that $r$ is at least as big as the largest power of a variable occurring
in a minimal generator of $J$.  However, now we choose the entries of $A$
generically in $R = K[X_1,\dots,X_n,u_1,\dots,u_t]$, where we even allow $t =
0$.  The degree of genericity we require is that the polynomials $F_j =
\prod_{i=1}^N L_{j,i}$, $1 \leq j \leq n$, define a complete intersection,
$X$.  Note that $F_j$ is the product of the entries of the $j$-th row, and
that the height of the complete intersection is $n$, the number of variables
in $S$.

The same construction as before produces from $J$ and $A$ an ideal $I$ of
$R$.  Now $I$ will no longer be a lifting in the sense of Definition
\ref{lift-def}.  However, we claim that Proposition \ref{same-resolution}
still holds, and hence so do any of the results of this section that do not
have to do with lifting.  As indicated in the proof of Proposition
\ref{same-resolution}, the first step is to observe that in any case Taylor's
resolution leads to a complex when we replace $J$ by $I$.  We again just have
to show that the Buchsbaum-Eisenbud exactness criterion gives that we have a
resolution.  More precisely, we have to show that 
\begin{itemize}
\item[(i)] $\rank \bar d_{k+1} + \rank \bar d_k = \rank \bar F_k$ for all $k$,
and 
\item[(ii)] the ideal of $(\rank \bar d_k )$-minors of $\bar d_k$ contains a
regular sequence of length $k$, or is equal to $R$.
\end{itemize}
Now, however, we do not
have that the restriction of $\bar d_k$ is $d_k$.  Nevertheless, the matrices
$\bar d_k$ still have the same form as the $d_k$.  The main thing to check is
that $\rank d_k = \rank \bar d_k$.

The rank of such a matrix
$d_k$ is the largest number of linearly independent columns, so clearly
$\rank \bar d_k \geq \rank d_k$.  For the reverse inequality, the danger is
that a linear combination of columns of $d_k$ could be zero while the
corresponding linear combinations of columns of $\bar d_k$ be non-zero.  This
could happen if $0 \neq c_{AB} = c_{A'B}$, where $B = \{i_1,\dots,i_{s-1} \}$,
$A = B \cup \{i_k \}$ and $A' = B \cup \{ i_{k}' \}$, but the analogous
entries for $I$ are not equal.  
Suppose
\[
c_{AB}  =  (-1)^k \frac{\lcm \{ m_{i_1} ,\dots,m_{i_{s-1}} ,m_{i_k} \}}
{\lcm \{ m_{i_1},\dots,m_{i_{s-1}} \} } \ \ \ \hbox{and} \ \ \ 
c_{A'B}  =  (-1)^k \frac{\lcm \{ m_{i_1} ,\dots,m_{i_{s-1}} ,m_{i_k'} \}}
{\lcm \{ m_{i_1},\dots,m_{i_{s-1}} \} }
\]
If $m_{i_k}$ contributes to the $\lcm$ then it contains a power of a variable
$X_j$ which is larger than the powers of $X_j$ in the other monomials
$m_{i_1},\dots,m_{i_{s-1}}$.  Since $c_{AB} = c_{A'B}$, $m_{i_k}'$ contains
the same power of $X_j$.  Hence they contribute the same number of entries
from the $j$-th row of the ``lifting'' matrix $A$, so the entries of the
relevant columns of $\bar d_k$ are equal.

From this fact, condition (i) follows immediately.  Condition (ii) follows
because any $(\rank d_k)$-minor of $d_k$ corresponds to a $(\rank \bar
d_k)$-minor of $\bar d_k$, and if $k$ of the former form a regular sequence
then clearly $k$ of the latter do as well.  
\end{remark}

Thanks to Remark \ref{more-general}, we now extend the notion of lifting to a
more general one.  See also Theorem \ref{lift-artinian}  and Corollary
\ref{lifting-case}.

\begin{definition} \label{pseudo-lifting-def}
Let $J$ be a monomial ideal in 
\[
S = K[X_1,\dots,X_n] \subset R =
k[X_1,\dots,X_n,u_1,\dots,u_t],
\]
 where $t \geq 0$.  Let $A$ be an $n
\times r$ matrix of linear forms in $R$, where $r$ is at least as big as the
largest power of a variable occurring in a minimal generator of $J$, and such
that the entries of $A$ satisfy the condition that the polynomials $F_j =
\prod_{i=1}^N L_{j,i}$, $1 \leq j \leq n$, define a complete intersection,
$X$. Let $I$ be the ideal obtained from $J$ by the construction described in
this section.  Then we shall call $I$ a {\em pseudo-lifting} of $J$.
\end{definition}


\section{Configurations of Linear Varieties}

The results of the preceding section give a number of nice properties of the
ideal $I$ obtained from the monomial ideal $J$, with no assumption on the
lifting matrix (or on $J$ other than being a monomial ideal).  In particular,
Corollary \ref{some-facts} says that
$I$ is a $t$-lifting of $J$.   We need one more fact in order to show that $I$
is a {\em reduced} $t$-lifting of $J$, namely that $I$ is radical.  To get
this, we need to make some assumptions on the lifting matrix.  We will prove
something more.  We will prove not only that $I$ defines the union, $V$, of
linear varieties, but in fact we would like to control the way that these
linear varieties intersect, and show that we can arrange that they intersect
in a very nice way.  For curves the ideal result would be to produce
so-called {\em stick figures}, i.e.\ unions of lines such that no more than
two pass through any given point.  

Consider however the following example, which shows that stick figures are too
ambitious in general, without some assumption on $J$, and also shows the
approach we will take.

\begin{example} \label{not-stick}
Let $n=3$ and consider the ideal 
\[
J = (X_1^3, X_1^2X_2, X_1^2X_3, X_1X_2^2, X_1X_2X_3, X_2^3, X_2^2X_3).
\] 
Note that $J$ is not saturated (see below) and that the saturation of $J$ is
not radical, defining instead a zeroscheme of degree 3 in $\proj{2}$
supported on a point $P$.  Let $t=1$. For the lifting, we will have to make
some ``generality'' assumption on the lifting matrix if we want to get a
radical ideal.  For example, if we took $L_{j,i} = X_j$ for all $1 \leq j
\leq 3$ and $1 \leq i \leq 3$, we see that $V$ is just a cone over the scheme
$W$ defined by $J$, hence not reduced.

Hence we will now assume that the linear forms $L_{j,i} \in
K[X_j,u_1]$ ($1 \leq j \leq n$) are chosen generally.  Then lifting we will
check that we obtain the saturated ideal of the union, $V$, of three lines in
$\proj{3}$ passing through $P$, together with three distinct points in
$\proj{3}$.  In particular, the top dimensional part is not a stick figure.

Note that $J = (X_1,X_2)^2 \cap (X_1,X_2,X_3)^3$.  Considering only the
ideal $(X_1,X_2)^2$, we lift to the saturated ideal of just the union of the
three lines passing through $P$.  If we consider instead the lifting of
$(X_1,X_2,X_3)^3$, we obtain instead the saturated ideal of 10 points in
$\proj{3}$.  However, 7 of these points lie on the union of the three lines,
none at the vertex, so that taking the union (i.e.\ intersecting the ideals)
gives the union of the three lines and three points described above.  (Since
7 is not divisible by 3, this fact is somewhat surprising: it says that the
three lines are not indistinguishable.)

To check this, following Remark \ref{primary-decomp-of-mon-id} and removing
redundant terms, notice that 
\[
(X_1,X_2)^2 = (X_1,X_2^2) \cap (X_1^2,X_2)
\]
 and 
\[
\begin{array}{rcl}
(X_1,X_2,X_3)^3
& =  & (X_1,X_2^2,X_3^2) \cap  (X_1^3,X_2,X_3) \cap (X_1,X_2,X_3^3) \cap
(X_1,X_2^3,X_3) \cap \\
&&  (X_1^2,X_2^2,X_3) \cap (X_1^2,X_2,X_3^2).
\end{array}
\]
In this
form it is easy to see what the components of $V$ will be after we lift:
again we remove redundant terms and we obtain that $(X_1,X_2)^2$ lifts to
\[
(L_{1,1},L_{2,1}) \cap (L_{1,1},L_{2,2}) \cap (L_{1,2},L_{2,1})
\]
while $(X_1,X_2,X_3)^3$ lifts to
\[
\begin{array}{c}
(L_{1,1},L_{2,1},L_{3,1}) \cap

(L_{1,1},L_{2,1},L_{3,2}) \cap

(L_{1,1},L_{2,2},L_{3,1}) \cap

(L_{1,1},L_{2,2},L_{3,2}) \cap \\

(L_{1,2},L_{2,1},L_{3,1}) \cap

(L_{1,3},L_{2,1},L_{3,1}) \cap

(L_{1,1},L_{2,1},L_{3,3}) \cap

(L_{1,1},L_{2,3},L_{3,1}) \cap \\

(L_{1,2},L_{2,2},L_{3,1}) \cap

(L_{1,2},L_{2,1},L_{3,2}).
\end{array}
\]
Then $J$ lifts to
\[
\begin{array}{c}
(L_{1,1},L_{2,1}) \cap (L_{1,1},L_{2,2}) \cap (L_{1,2},L_{2,1}) \\
\cap \\
(L_{1,3},L_{2,1},L_{3,1}) \cap

(L_{1,1},L_{2,3},L_{3,1}) \cap 

(L_{1,2},L_{2,2},L_{3,1}) 
\end{array}
\]
as claimed.
\end{example}

Example \ref{not-stick} shows that without some assumptions on $J$ we cannot
hope to get stick figures, but that we can obtain reduced unions of linear
varieties.  For the remainder of this paper we make the following convention
(but see also Remark \ref{matrix-condition}):

\begin{quotation}
{\em For the lifting matrix $A = [L_{j,i}]$, we assume that on each row
(corresponding to a variable $X_j$), each entry $L_{j,i}$ is chosen
generically with respect to the preceding entries $L_{j,1},\dots,L_{j,i-1}$
on that row.}
\end{quotation}
Note, however, that for this section we continue to assume that $L_{j,i} \in
K[X_j,u_1,\dots,u_t]$.  In the next section we will use more generally chosen
linear forms, thanks to Remark \ref{more-general}.

\begin{remark} \label{matrix-condition}
We would like to make somewhat more precise what we mean by ``generically''
in the above assumption.  Geometrically, the condition is as follows:  For
any choice of $k$ pairwise distinct entries $L_1,\dots,L_k$ from the matrix,
we require
\begin{equation} \label{codim-cond}
\codim (L_1,\dots,L_k) = \min \{ n+t,k \}.
\end{equation}
For any given monomial ideal $J$, note that this puts only a finite number of
conditions on the lifting matrix.  In fact, let $N_k$ ($1 \leq k \leq n$) be
the largest power of $X_k$ occurring as a factor of a minimal generator of
$J$.  Let $N = \max \{ N_1,\dots,N_n \}$.  Then the lifting matrix $A=
[a_{k,\ell}]$ for $J$ can be chosen to have size $n \times N$, and we can
further assume that
\[
a_{k,\ell} = 0 \ \ \hbox{ if } \ell > N_k.
\]
Then the above matrix condition is only required for non-zero entries of $A$.

We can describe a suitable matrix $A$ explicitly.  Put $p = N_1+\dots+ N_n$
and choose $p$ distinct elements $b_1,\dots,b_p \in K$.  Consider the
Vandermonde matrix
\[
B = \left [
\begin{array}{cccccccc}
1 & 1 & \cdots & 1 \\
b_1 & b_2 & \cdots & b_p \\
\vdots & \vdots & \cdots & \vdots \\
b_1^t & b_2^t & \cdots & b_p^t
\end{array}
\right ]
\]
Now let
\[
\begin{array}{rcl}
B_1 & = & \hbox{submatrix formed by the first $N_1$ columns of $B$} \\
B_2 & = & \hbox{submatrix formed by the next $N_2$ columns of $B$} \\
& \vdots \\
B_n & = & \hbox{submatrix formed by the next $N_n$ columns of $B$} 
\end{array}
\]
Now we produce a lifting matrix $A$, of size $n \times N$, as follows.  The
first $N_k$ entries of the $k$-th row are given by
\[
(x_k,u_1,\dots,u_t) \cdot B_k.
\]
If $N_k < N$, we ``complete'' the $k$-th row by zeros.  Then the properties
of the Vandermonde matrix ensure that the matrix condition (\ref{codim-cond})
holds true for $A$ if we choose only non-zero entries of $A$.  It follows
that we can lift the monomial ideal $J$ to a ``nice'' ideal $I$ provided our
ground field $K$ has at least $p$ elements.

We have already seen in Example \ref{not-stick}, and it will be made more
precise in this section and especially in section 4, that the components of
the scheme defined by the lifted ideal are linear varieties defined by
suitable subsets of the entries of the lifting matrix (at most one entry from
each row).  Let $\wp_1,\dots,\wp_j$ be associated prime ideals of the lifted
ideal $I$ of $J$ using a lifting matrix $A$.  Then it follows from
(\ref{codim-cond}) that 
\[
\begin{array}{rcl}
\codim (\wp_1 + \dots + \wp_j) & = & \min \{ n+t, \hbox{\# entries of the
lifting matrix occurring} \\
&& \phantom{\min \{ n+t,} \hbox{as minimal generators of some $\wp_i$ } \}
\end{array}
\]
The same holds for the explicit lifting matrix given above.  Hence all the
results below about generalized stick figures hold when the lifting matrix is
this explicit matrix.

In the case $t=1$ and $J$ Artinian, this matrix was essentially given in
\cite{GGR}.
\end{remark}

Inspired by Example \ref{not-stick}, we now would like to examine the lifting
more carefully.  We would like to show, first, that we have a reduced
$t$-lifting (cf.\ Definition \ref{t-lifting-def}).  But in fact, we would
like to see when we have a stick figure, for curves, and to extend this
notion to higher dimension.  In \cite{BM5} the authors defined a ``good linear
configuration'' to be a locally Cohen-Macaulay codimension two union of
linear subspaces in $\proj{n}$ such that the intersection of any three
components has dimension at most $n-4$.  We would like to modify that
definition slightly (removing the locally Cohen-Macaulay assumption and
allowing arbitrary codimension), as follows.  See also Remark
\ref{unambitious}.

\begin{definition} \label{gen-stick-fig}
Let $V$ be a union of linear subvarieties of $\proj{m}$ of the same dimension
$d$.  Then $V$ is a generalized stick figure if the intersection of any
three components of $V$ has dimension at most $d-2$ (where the empty set is
taken to have dimension $-1$).  In particular, if
$d=1$ then $V$ is a stick figure.
\end{definition}

\begin{theorem}  \label{lift-to-stick}
Let $J \subset S$ be a monomial ideal, let $I$ be the $t$-lifting of $J$ 
described in the paragraph preceding Remark \ref{geometrically},
and let $V$ (resp.\ $W$) be the subscheme of
$\proj{n+t-1}$ defined by $I$ (resp.\ the subscheme of $\proj{n-1}$ defined by
$J$).  Then 
\begin{itemize}
\item[(a)] $I$ is a radical ideal in $R$, i.e.\ $I$ is a reduced $t$-lifting
of $J$; in fact,  $V$ is a union of linear varieties.  
\item[(b)] The union of the components of $V$ of any given dimension $d$ form
a generalized stick figure away from $W$, in the sense that 
\[
(\hbox{the intersection of any 3 components}) \backslash W_{red}
\]
has dimension at most $d-2$.

\item[(c)] Suppose that $J$ is Artinian.  Then $V$ is an arithmetically 
Cohen-Macaulay generalized stick figure.

\item[(d)] Suppose that $J$ is unmixed and not Artinian, and let $J = Q_1 \cap
\dots \cap Q_r$ be a minimal primary decomposition of $J$.  Let $W_i$ be the
scheme defined by $Q_i$, for each $i$.  If $t=1$ then $V$ is a generalized
stick figure if and only if $\deg W_i \leq 2$ for all $i$.  If $t \geq 2$
then $V$ is a generalized stick figure.  
\end{itemize}
\end{theorem}

\begin{proof}

We first prove that $V$ is a union of linear varieties, and hence that $I$ is
radical.  We know from Remark \ref{primary-decomp-of-mon-id} that we can
write $J$ in the form
$J = Q_1 \cap \dots \cap Q_r$ where each
$Q_i$ is a complete intersection of the form $(X_{i_1}^{a_1} ,\dots,
X_{i_p}^{a_p} )$.  Using Remark \ref{geometrically}, it is clear that each
lifted ideal $\bar Q_i$ defines a reduced complete intersection in
$\proj{n+t-1}$.  Finally, from Lemma \ref{decompose-I} we get that $I$ is
reduced and a union of linear varieties.

We have seen in Example \ref{not-stick} that we cannot hope that $V$ will
always be a generalized stick figure.  Indeed, by Remark \ref{geometrically},
if $W$ is not empty then all components of $V$ will pass through $W$, so if
there are too many of these components then it will fail to be a generalized
stick figure.  

For (b), consider the lifting matrix 
\[
A = 
\left [
\begin{array}{ccccccccc}
L_{1,1} & L_{1,2} & L_{1,3} & \dots \\
L_{2,1} & L_{2,2} & L_{2,3} & \dots \\
\vdots & \vdots & \vdots \\
L_{n,1} & L_{n,2} & L_{n,3} & \dots 
\end{array}
\right ].
\]
Consider a primary decomposition of $J$ and write it as $J = J_1 \cap J_2
\cap \dots \cap J_s$ where, for $1 \leq c \leq p$, $J_c$ is the intersection
of the primary components of codimension $c$.  By Lemma \ref{decompose-I} and
Corollary \ref{unmixed}, the components of $V$ of codimension $c$ are
defined by the lifting $I_c$ of $J_c$.  Let $W_c$ be the scheme defined by
$J_c$ (possibly empty if $J_c$ is Artinian) and let $V_c$ be the scheme
defined by $I_c$.  

Any component of $V_c$ is a linear variety defined by $c$ entries of the
lifting matrix, no two from any given row, and it vanishes (set
theoretically) on a component of $W_c$.  Because the linear forms in $A$ were
chosen generally, subject only to the condition that the entries of the
$j$-th row are linear forms in the variables $X_j,u_1,\dots,u_t$, it follows
that the intersection of any $3$ components has codimension $2$ in $V_c$,
at least away from the linear space $\proj{n-1}$ defined by the vanishing of
the variables $u_1,\dots,u_n$.  (See also Remark \ref{unambitious}.)  But the
intersection of $V$ with this linear space is exactly $W$.  This proves (b).

For (c), the fact that $V$ is \acm follows from Corollary \ref{some-facts}, 
while the fact that it is a generalized stick figure follows from (b), since
$W$ is empty.

For (d), let $V_i$ be the lifting of $W_i$.  Note that $\deg V_i = \deg W_i$,
by Corollary \ref{some-facts} (iv), $V$ is unmixed by Corollary
\ref{unmixed}, and hence $V$ has exactly $\deg W_i$ components passing
through $W_i$.  Furthermore, $W_i$ has codimension $t$ in $V$.  Then using
(b), we see that the only way that the intersection of three components of
$V$ could fail to have codimension two in $V$ is if $t=1$ and $\deg W_i \geq
3$ for some $i$.
\end{proof}

\begin{remark}
If $J$ is Artinian, the fact that $I$ is radical follows from \cite{GGR}.
Indeed, they show that a 1-lifting is a reduced set of points, and it is
clear from the above that in the case of a $t$-lifting, a sequence of
hyperplane sections reduces to a 1-lifting.  If this 1-lifting is reduced
then the original $t$-lifting must be reduced.
\end{remark}

\begin{remark}\label{unambitious}
One might wonder why, in Definition \ref{gen-stick-fig}, we do not define a
generalized stick figure to involve the intersection of more than three
components.  As a first answer, let $R$ be the polynomial ring in 5 variables
and let $A_1,A_2,A_3, A_4$ be general linear forms.  Let $F = A_1\cdot A_2$
and let $G = A_3\cdot A_4$.  The complete intersection of $F$ and $G$ defines
a union of four linear varieties of dimension 2.  Any reasonable definition
of a generalized stick figure must, in our opinion, include this example. 
However, the intersection of these four components has dimension 0 rather
than being empty.

More generally, in the context of lifting, consider a lifting matrix
\[
A = 
\left [
\begin{array}{ccccccc}
L_{1,1} & L_{1,2} & \dots & L_{1,r} \\
L_{2,1} & L_{2,2} & \dots & L_{2,r} \\
\vdots & \vdots &  &\vdots \\
L_{n,1} & L_{n,2} & \dots & L_{n,r} 
\end{array}
\right ]
\]
and let $J = (X_1^{a_1},X_2^{a_2},X_3^{a_2},\dots,X_n^{a_n})$ with $a_1 \geq
2$ and $a_2 \geq 2$. $J$ is Artinian, and the lifted ideal $I$ defines a
codimension $n$ union of varieties, $V$, which contains in particular the
components $V_i$, $1 \leq i \leq 4$, defined by the ideals 
\[
\begin{array}{ccccccc}
I_{V_1} = (L_{1,1},L_{2,1}, L_{3,1},\dots,L_{n,1}) \\
I_{V_2} = (L_{1,1},L_{2,2}, L_{3,1},\dots,L_{n,1}) \\
I_{V_3} = (L_{1,2},L_{2,1}, L_{3,1},\dots,L_{n,1}) \\
I_{V_4} = (L_{1,2},L_{2,2}, L_{3,1},\dots,L_{n,1}) .
\end{array}
\]
Assume that $t \geq 3$, so that $\dim V = t-1 \geq 2$ in $\proj{n+t-1}$.  But
then $\dim V_1 \cap V_2 \cap V_3 \cap V_4 = t-3$, not $t-4$.  Hence for
part (c) of Theorem \ref{lift-to-stick} to be true, we must avoid
intersecting four components.
\end{remark}

For the last result of this section, we recall some basic facts about
so-called ``O-sequences,'' which can be found in \cite{GMR}.  Macaulay showed
that a sequence of non-negative integers $\{ c_i \}, \ i \geq 0$, can be the
Hilbert function of a graded algebra $A = S/I$ if and only if $c_0 = 1$ and
the $c_i$ satisfy a certain growth condition (cf.\ \cite{stanley}).  Such a
sequence is called an {\em O-sequence}.  The sequence is said to be {\em
differentiable} if the first difference sequence $\{b_i \} = \Delta \{ c_i
\}$, defined by $b_i = c_i - c_{i-1}$, is also an O-sequence.  (We adopt the
convention that $c_{-1} = 0$.)  Extending this gives the notion of a {\em
$t$-times differentiable O-sequence}, in case $t$ successive differences
$\Delta^j \{ c_i \}$, $1 \leq j \leq t$, are all  still O-sequences.  In
\cite{GMR} the authors showed that any differentiable O-sequence occurs as
the Hilbert function of some graded algebra $S/I$, where
$I$ is {\em radical}.

An O-sequence $\{ c_i\}$ is said to have {\em dimension $d \geq 1$} (here
``dimension'' should be thought of as Krull dimension) if there  is a non-zero
polynomial $f(x)$, with rational coefficients, of degree $d-1$, such that for
all $s \gg 0$, $f(s) = c_s$.  If $c_i = 0$ for all $i \gg 0$ then we say
that $\{c_i\}$ has {\em dimension 0}.  The sequence $\{c_i\}$ is a
$d$-times differentiable O-sequence of dimension $d$ if and only if it is the
Hilbert function of a Cohen-Macaulay algebra $A = S/I$ of Krull dimension
$d$.  In this case the $d$-th successive difference $\Delta^d \{ c_i \}$ is
the Hilbert function of the {\em Artinian reduction} of $A$.  It is a finite
sequence of positive integers, called the {\em $h$-vector} of $A$ (cf.\ for
instance \cite{migliore}).  By \cite{GMR}, $I$ can be taken to be radical.  We
extend this with the following result, which was known in codimension two.

\begin{corollary}\label{differentiable}
Let $\{ c_i \}$ be a $t$-times differentiable O-sequence of dimension $t$. 
Then $\{c_i \}$ is the Hilbert function of an arithmetically Cohen-Macaulay
generalized stick figure of dimension $t-1$ (as a subscheme of projective
space).  In particular, any 2-differentiable O-sequence of dimension 2 is the
Hilbert function of some stick figure curve in projective space.
\end{corollary}

Inverse Gr\"obner basis theory asks which monomial ideals arise as initial
ideals of classes of ideals with prescribed properties.  For example, it is
conjectured that certain monomial ideals cannot be the initial ideal of a
prime ideal.  Instead, if we allow radical ideals, we have the following
result.

\begin{corollary}
Let $J \subset S$ be a monomial ideal.  Let $>$ be the degree-lex order such
that $X_1 > \dots > X_n > u_1 > \dots > u_t$.  Then there is a radical ideal
$I$ in $R$ such that the initial ideal of $I$ with respect to $>$ is $J \cdot
R$.  Moreover, if $J$ is an  Artinian monomial ideal in $S$ then $I$ can be
chosen to be the defining ideal of an \acm generalized stick figure.
\end{corollary}

\begin{proof}
Let $I \subset R$ be a reduced $t$-lifting of $J$.  Then we clearly have $J
\cdot R
\subset \hbox{in}(I)$.  On the other hand we have for the Hilbert functions: 
$\Delta^t h_{R/I}(j) = h_{S/J}(j) = \Delta^t h_{R/J \cdot R}(j)$.  Thus, the
equality $J \cdot R = \hbox{in}(I)$ follows. 
\end{proof}


\section{Pseudo-Liftings and Liftings of Artinian monomial ideals}

One of the main goals of this section is to describe the configurations of
linear varieties which arise by lifting Artinian monomial ideals, in the
sense of Definition \ref{t-lifting-def}.  We give the answer in  Corollary
\ref{lifting-case}.  This is actually a special case of the more general
notion of pseudo-lifting, as indicated in Remark \ref{more-general}, and we
describe this situation in Theorem \ref{lift-artinian}.

In this section we will consider a matrix 
\[
A = 
\left [
\begin{array}{ccccccccc}
L_{1,1} & L_{1,2} & L_{1,3} & \dots & L_{1,N} \\
L_{2,1} & L_{2,2} & L_{2,3} & \dots & L_{2,N} \\
\vdots & \vdots & \vdots & & \vdots \\
L_{n,1} & L_{n,2} & L_{n,3} & \dots & L_{n,N}
\end{array}
\right ]
\]
of linear forms in $R = K[X_1,\dots,X_n,u_1,\dots,u_t]$ satisfying
the following genericity property, slightly more than we assumed in Remark
\ref{more-general}:

\begin{quotation}
Let $F_j$ be the product of the entries of the $j$-th row of $A$.  We
assume that the ideal $(F_1,\dots,F_n)$ defines a reduced complete
intersection.
\end{quotation}
 
Throughout this section, ``configuration'' shall mean a reduced finite union
of linear varieties all of the same dimension.   We would like to describe
geometrically the configurations which are the pseudo-liftings of Artinian
monomial ideals, in the sense of Remark \ref{more-general} and Definition
\ref{pseudo-lifting-def}.  Recall that if the  entries $L_{j,i}$ are not in
$K[X_j,u_1,\dots,u_n]$,  these are not true liftings in the sense of
Definition \ref{t-lifting-def}.  

We have already seen in Theorem \ref{lift-to-stick} (c) that these
configurations will be \acm generalized stick figures.  In the case of
zero-dimensional schemes, the similarity of our construction to the notion of
a $k$-configuration will be evident (but they are not quite the same), and we
discuss the relation in Remark \ref{k-config}.  Part of the discussion will
involve the special case of Artinian lex-segment ideals, which we now recall.

\begin{definition}
Let $>$ denote the degree-lexicographic order on monomial ideals, i.e.\
$x_1^{a_1}\cdots x_n^{a_n} > x_1^{b_1}\cdots x_n^{b_n}$ if the first nonzero
coordinate of the vector
\[
\left ( \sum_{i=1}^n (a_i - b_i), a_1 - b_1 ,\dots,a_n - b_n \right )
\]
is positive.  Let $J$ be a monomial ideal.  Let $m_1,m_2$ be monomials in
$S$ of the same degree such that $m_1 > m_2$.  Then $J$ is a lex-segment
ideal if $m_2 \in J$ implies $m_1 \in J$.
\end{definition}

Let $J$ be an Artinian monomial ideal.  Let $N_j$ be the maximum power of
$X_j$ that occurs in a minimal generator of $J$, and let $N$ be the maximum
of the $N_j$.  Consider the  matrix
\[
A = 
\left [
\begin{array}{ccccccccc}
L_{1,1} & L_{1,2} & L_{1,3} & \dots & L_{1,N} \\
L_{2,1} & L_{2,2} & L_{2,3} & \dots & L_{2,N} \\
\vdots & \vdots & \vdots & & \vdots \\
L_{n,1} & L_{n,2} & L_{n,3} & \dots & L_{n,N}
\end{array}
\right ]
\]
where $L_{j,i}$ is a generally chosen linear form in the ring $R =
K[X_1,\dots,X_n,u_1,\dots,u_t]$.  Let
$I$ be the ideal obtained from $J$ using $A$, as in Remark \ref{more-general},
and let $V$ be the configuration of linear varieties defined by the saturated
ideal $I$.  Note that the dimension of $V$ as a projective subscheme is
$t-1$.  Let $F_1,\dots,F_n$ be defined by $F_j = \prod_{i=1}^N L_{j,i}$. 
Clearly $F_j \in I$ for all $j$.  

On the geometric side, the complete intersection $(F_1,\dots,F_n)$ defines a
union, $X$, of linear varieties of dimension $t-1$ and $V$ is a subset of
$X$.  Each component of $X$ is given by an $n$-tuple of linear forms
$(L_{1,i_1},\dots,L_{n,i_n})$.  Let $\Lambda \subset X$ be the component
corresponding to $(L_{1,i_1},\dots,L_{n,i_n})$.  By abuse of notation we will
write $\Lambda = (L_{1,i_1},\dots,L_{n,i_n}) \in V$ to mean that the
corresponding linear variety is a component of $V$.  Then recalling from
(\ref{lift-m-def}) how the pseudo-lifting is defined, we see that 
\begin{equation}\label{condition1}
\begin{array}{rcl}
(L_{1,i_1},\dots,L_{n,i_n}) \in V & \Leftrightarrow & \hbox{every monomial
generator $m \in J$ is divisible by} \\ 
&& \hbox{at least one of $X_1^{i_1} ,X_2^{i_2}, \dots ,X_n^{i_n}$.}
\end{array}
\end{equation}

As a result of (\ref{condition1}) we immediately get the following
description of the components of $V$, which is essentially contained in
\cite{GGR} in the case of dimension zero.

\begin{lemma}\label{not-in-J}
$(L_{1,i_1},\dots,L_{n,i_n}) \in V$ if and only if $X_1^{i_1-1}\cdots
X_n^{i_n-1} \notin J.$
\end{lemma}

We would like to give a geometric description of the configurations that
arise in this way.  

Since $J$ is an Artinian monomial ideal, it contains powers of all the
variables $X_j$.  Let $\alpha_j$ be the least degree of $X_j$ that
appears in $J$.  Without loss of generality, order and relabel the $X_j$ so
that $\alpha_1 \leq \alpha_2 \leq \dots \leq \alpha_n$.  

Notice that, by (\ref{condition1}), for $i \geq 2$ we have that if the
component 
$(L_{1,i},L_{2,i_2},\dots,L_{n,i_n})$ is in $V$ then the component 
$(L_{1,i-1},L_{2,i_2},\dots,L_{n,i_n})$ is in $V$.  More generally, we
immediately see the following:

\begin{lemma}\label{condition2}
For $1 \leq j \leq n$ and $i_j \geq 2$, if $(L_{1,i_1},L_{2,i_2},\dots,
L_{j,i_j}, \dots, L_{n,i_n}) \in V$ then $(L_{1,i_1},L_{2,i_2},\dots,
L_{j,i_j -1}, \dots, L_{n,i_n}) \in V$.
\end{lemma}

Furthermore, in the special case where $J$ is a lex-segment ideal we have a
stronger property:

\begin{lemma}\label{condition3}
Let $J$ be an Artinian monomial lex-segment ideal.  For $1 \leq j \leq n$ and
$i_j \geq 2$, if $(L_{1,i_1},L_{2,i_2},\dots, L_{j,i_j}, \dots, L_{n,i_n}) \in
V$ then 
\[
(L_{1,i_1},L_{2,i_2},\dots, L_{j,i_j -1},  L_{j+1,i_{j+1}'}
,\dots,L_{n,i_n'} ) \in V 
\]
for any $i_{j+1}',\dots,i_n' \in {\mathbb N}$ such that $i_{j+1}' +
\cdots+i_n' = i_{j+1} + \cdots + i_n + 1$.
\end{lemma}

\begin{proof}
Since $(L_{1,i_1},L_{2,i_2},\dots, L_{j,i_j}, \dots, L_{n,i_n}) \in V$, we
have by Lemma \ref{not-in-J} that 
\[
X_1^{i_1 -1} \cdots X_j^{i_j-1} \cdots X_n^{i_n -1} \notin J.
\]
But $X_1^{i_1 -1} \cdots X_j^{i_j-1} \cdots X_n^{i_n -1} > 
X_1^{i_1 -1} \cdots X_j^{i_j-2} X_{j+1}^{i_{j+1}'-1} \cdots X_n^{i_n' -1}$,
so  since $J$ is a lex-segment ideal we get that 
\[
X_1^{i_1 -1} \cdots X_j^{i_j-2} X_{j+1}^{i_{j+1}'-1} \cdots X_n^{i_n' -1}
\notin J
\]
as well.  Hence 
\[
(L_{1,i_1},L_{2,i_2},\dots, L_{j,i_j -1},  L_{j+1,i_{j+1}'}
,\dots,L_{n,i_n'} ) \in V
\]
as claimed.
\end{proof}

\begin{example} \label{useful-ex}

Let 
\[
J = (X_1^3,X_1^2X_2^2, X_1^2X_2X_3, X_1 X_2^3, X_1X_2^2X_3, X_1X_2X_3^2,
X_1X_3^3, X_2^4, X_2^3X_3, X_2^2X_3^2, X_2X_3^3 X_3^4 ).
\]
Note that $J$ is ``almost'' a lex-segment ideal, missing only the monomial
$X_1^2X_3^2$.  The $h$-vector of $J$ is $(1,3,6,9,1)$.  We first sketch the
configuration obtained by pseudo-lifting $J$, by showing the configurations on
the $L_{1,1}$-plane, the $L_{1,2}$-plane and the $L_{1,3}$-plane.

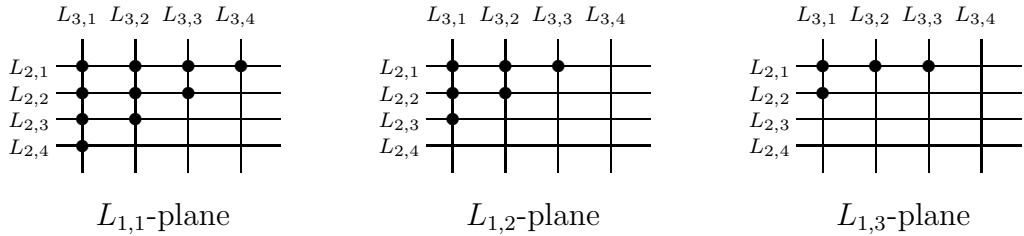
\newsavebox{\running}
\savebox{\running}(100,85)[tl]
{
\put (120,30){\line (0,1){50}}
\put (140,30){\line (0,1){50}}
\put (160,30){\line (0,1){50}}
\put (180,30){\line (0,1){50}}

\put (110,40){\line (1,0){85}}
\put (110,50){\line (1,0){85}}
\put (110,60){\line (1,0){85}}
\put (110,70){\line (1,0){85}}

\put (117,37){$\bullet$}
\put (117,47){$\bullet$}
\put (117,57){$\bullet$}
\put (117,67){$\bullet$}
\put (137,67){$\bullet$}
\put (137,47){$\bullet$}
\put (137,57){$\bullet$}
\put (157,57){$\bullet$}
\put (157,67){$\bullet$}
\put (177,67){$\bullet$}

\put (92,37){$\scriptstyle L_{2,4}$}
\put (92,47){$\scriptstyle L_{2,3}$}
\put (92,57){$\scriptstyle L_{2,2}$}
\put (92,67){$\scriptstyle L_{2,1}$}

\put (110,87){$\scriptstyle L_{3,1}$}
\put (130,87){$\scriptstyle L_{3,2}$}
\put (150,87){$\scriptstyle L_{3,3}$}
\put (170,87){$\scriptstyle L_{3,4}$}

\put (125,10){$L_{1,1}$-plane}
}

\newsavebox{\runningg}
\savebox{\runningg}(100,85)[tl]
{
\put (120,30){\line (0,1){50}}
\put (140,30){\line (0,1){50}}
\put (160,30){\line (0,1){50}}
\put (180,30){\line (0,1){50}}

\put (110,40){\line (1,0){85}}
\put (110,50){\line (1,0){85}}
\put (110,60){\line (1,0){85}}
\put (110,70){\line (1,0){85}}

\put (117,47){$\bullet$}
\put (117,57){$\bullet$}
\put (117,67){$\bullet$}
\put (137,67){$\bullet$}
\put (137,57){$\bullet$}
\put (157,67){$\bullet$}

\put (92,37){$\scriptstyle L_{2,4}$}
\put (92,47){$\scriptstyle L_{2,3}$}
\put (92,57){$\scriptstyle L_{2,2}$}
\put (92,67){$\scriptstyle L_{2,1}$}

\put (110,87){$\scriptstyle L_{3,1}$}
\put (130,87){$\scriptstyle L_{3,2}$}
\put (150,87){$\scriptstyle L_{3,3}$}
\put (170,87){$\scriptstyle L_{3,4}$}

\put (125,10){$L_{1,2}$-plane}
}

\newsavebox{\runninggg}
\savebox{\runninggg}(100,85)[tl]
{
\put (120,30){\line (0,1){50}}
\put (140,30){\line (0,1){50}}
\put (160,30){\line (0,1){50}}
\put (180,30){\line (0,1){50}}

\put (110,40){\line (1,0){85}}
\put (110,50){\line (1,0){85}}
\put (110,60){\line (1,0){85}}
\put (110,70){\line (1,0){85}}

\put (117,57){$\bullet$}
\put (117,67){$\bullet$}
\put (137,67){$\bullet$}
\put (157,67){$\bullet$}

\put (92,37){$\scriptstyle L_{2,4}$}
\put (92,47){$\scriptstyle L_{2,3}$}
\put (92,57){$\scriptstyle L_{2,2}$}
\put (92,67){$\scriptstyle L_{2,1}$}

\put (110,87){$\scriptstyle L_{3,1}$}
\put (130,87){$\scriptstyle L_{3,2}$}
\put (150,87){$\scriptstyle L_{3,3}$}
\put (170,87){$\scriptstyle L_{3,4}$}

\put (125,10){$L_{1,3}$-plane}
}

\begin{picture}(400,105) 
\put(-80,0){\usebox{\running}}
\put(60,0){\usebox{\runningg}}
\put(200,0){\usebox{\runninggg}}

\end{picture}

\noindent This satisfies the condition of Lemma \ref{condition2}, of course,
but not Lemma \ref{condition3}. The ``offending'' point is
$(L_{1,3}, L_{2,1},L_{3,3})$,  since the point $(L_{1,2}, L_{2,1},L_{3,4})$
is not in the configuration (take $j=1$ and $k=3$).  Hence this configuration
does not occur for lex-segment ideals.
\end{example}

\begin{example}
Here is a slightly more subtle example.  Let $n=3$ and consider the monomial
ideal 
\[
J = (X_1^3, X_1^2X_2,X_1^2X_3,X_1X_2^2, X_1X_2
X_3,X_2^3,X_1X_3^3,X_2^2X_3^2,X_2X_3^3, X_3^4).
\]
The component in degree 3 fails to be lex-segment because it is missing the
monomial $X_1X_3^2$ after $X_1X_2X_3$.  The configuration corresponding to
the pseudo-lifting of $J$ looks as follows:

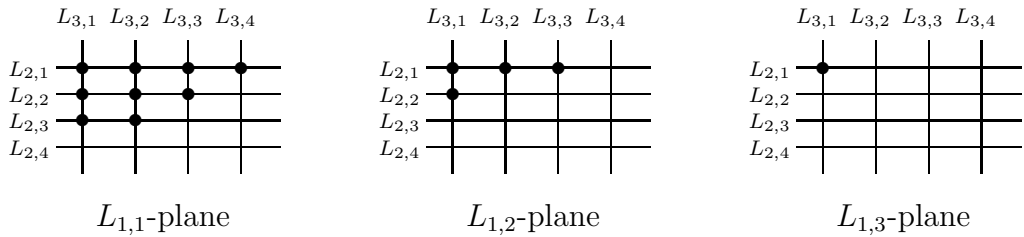
\newsavebox{\ctrex}
\savebox{\ctrex}(100,85)[tl]
{
\put (120,30){\line (0,1){50}}
\put (140,30){\line (0,1){50}}
\put (160,30){\line (0,1){50}}
\put (180,30){\line (0,1){50}}

\put (110,40){\line (1,0){85}}
\put (110,50){\line (1,0){85}}
\put (110,60){\line (1,0){85}}
\put (110,70){\line (1,0){85}}

\put (117,47){$\bullet$}
\put (117,57){$\bullet$}
\put (117,67){$\bullet$}
\put (137,67){$\bullet$}
\put (137,47){$\bullet$}
\put (137,57){$\bullet$}
\put (157,57){$\bullet$}
\put (157,67){$\bullet$}
\put (177,67){$\bullet$}

\put (92,37){$\scriptstyle L_{2,4}$}
\put (92,47){$\scriptstyle L_{2,3}$}
\put (92,57){$\scriptstyle L_{2,2}$}
\put (92,67){$\scriptstyle L_{2,1}$}

\put (110,87){$\scriptstyle L_{3,1}$}
\put (130,87){$\scriptstyle L_{3,2}$}
\put (150,87){$\scriptstyle L_{3,3}$}
\put (170,87){$\scriptstyle L_{3,4}$}

\put (125,10){$L_{1,1}$-plane}
}

\newsavebox{\ctrexx}
\savebox{\ctrexx}(100,85)[tl]
{
\put (120,30){\line (0,1){50}}
\put (140,30){\line (0,1){50}}
\put (160,30){\line (0,1){50}}
\put (180,30){\line (0,1){50}}

\put (110,40){\line (1,0){85}}
\put (110,50){\line (1,0){85}}
\put (110,60){\line (1,0){85}}
\put (110,70){\line (1,0){85}}

\put (117,57){$\bullet$}
\put (117,67){$\bullet$}
\put (137,67){$\bullet$}
\put (157,67){$\bullet$}

\put (92,37){$\scriptstyle L_{2,4}$}
\put (92,47){$\scriptstyle L_{2,3}$}
\put (92,57){$\scriptstyle L_{2,2}$}
\put (92,67){$\scriptstyle L_{2,1}$}

\put (110,87){$\scriptstyle L_{3,1}$}
\put (130,87){$\scriptstyle L_{3,2}$}
\put (150,87){$\scriptstyle L_{3,3}$}
\put (170,87){$\scriptstyle L_{3,4}$}

\put (125,10){$L_{1,2}$-plane}
}

\newsavebox{\ctrexxx}
\savebox{\ctrexxx}(100,85)[tl]
{
\put (120,30){\line (0,1){50}}
\put (140,30){\line (0,1){50}}
\put (160,30){\line (0,1){50}}
\put (180,30){\line (0,1){50}}

\put (110,40){\line (1,0){85}}
\put (110,50){\line (1,0){85}}
\put (110,60){\line (1,0){85}}
\put (110,70){\line (1,0){85}}

\put (117,67){$\bullet$}

\put (92,37){$\scriptstyle L_{2,4}$}
\put (92,47){$\scriptstyle L_{2,3}$}
\put (92,57){$\scriptstyle L_{2,2}$}
\put (92,67){$\scriptstyle L_{2,1}$}

\put (110,87){$\scriptstyle L_{3,1}$}
\put (130,87){$\scriptstyle L_{3,2}$}
\put (150,87){$\scriptstyle L_{3,3}$}
\put (170,87){$\scriptstyle L_{3,4}$}

\put (125,10){$L_{1,3}$-plane}
}

\begin{picture}(400,105) 
\put(-80,0){\usebox{\ctrex}}
\put(60,0){\usebox{\ctrexx}}
\put(200,0){\usebox{\ctrexxx}}

\end{picture}

At first glance one would be tempted to say that this satisfies the condition
of Lemma \ref{condition3}, just from an examination of the picture.  However,
the fact that the point \linebreak $(L_{1,2}, L_{2,1}, L_{3,3})$ is in $V$
requires that $(L_{1,1}, L_{2,4}, L_{3,1}) \in V$ as well, and this does not
hold. 
\end{example}

We can now describe the configurations which arise as pseudo-liftings of
Artinian monomial ideals and those which arise as pseudo-liftings of Artinian
lex-segment monomial ideals.

\begin{theorem}\label{lift-artinian}
Let $V$ be a configuration of linear varieties of codimension $n$ in
$\proj{n+t-1}$.   Let $A$ be a matrix of linear forms, 
\[
A = 
\left [
\begin{array}{ccccccccc}
L_{1,1} & L_{1,2} & L_{1,3} & \dots & L_{1,N} \\
L_{2,1} & L_{2,2} & L_{2,3} & \dots & L_{2,N} \\
\vdots & \vdots & \vdots & & \vdots \\
L_{n,1} & L_{n,2} & L_{n,3} & \dots & L_{n,N}
\end{array}
\right ]
\]
such that the polynomials $F_j = \prod_{i=1}^N L_{j,i}$, $1 \leq j \leq n$,
define a {\em reduced} complete intersection, $X$, containing $V$.
Then $V$ is the pseudo-lifting, via $A$, of an Artinian monomial ideal in $n$
variables if and only if the condition of Lemma \ref{condition2} holds. 
Furthermore, $V$ is the pseudo-lifting, via $A$, of an Artinian lex-segment
monomial ideal if and only if in addition the condition of Lemma
\ref{condition3} holds.
\end{theorem}

\begin{proof}
For both parts of the theorem we have only to prove the sufficiency of the
condition.  We continue to use the same notation $L_{j,i}$ for the linear
form and for the corresponding hyperplane.  Without loss of generality we can
assume that $V$ neither empty nor all of $X$.  

Let $I$ be the ideal generated
by the set of all polynomials which are products of the form
\[
\prod_{j=1}^n \left ( \prod_{i=1}^{a_j}  L_{j,i} \right )
\]
and which vanish on $V$.  Remove generators from this set to obtain a minimal
generating set.  Any generator of this form corresponds, via
(\ref{lift-m-def}), to a monomial.  Let $J$ be the ideal generated by the set
of all such monomials.  $J$ is Artinian since it contains a power of each of
the variables, and clearly $I$ is the pseudo-lifting of $J$.  Hence $I$ is the
saturated ideal (thanks to Corollary \ref{some-facts}) of a scheme $Z$
with $V \subseteq Z \subseteq X$. 

We have only to show that $Z = V$, and for this we use the condition of Lemma
\ref{condition2}.  Let $\Lambda$ be a component of $X$ which is not in $V$ (OK
since $V$ is not all of $X$).  We will produce a polynomial $F$ in $I$ which
does not vanish on $\Lambda$.  Being a component of $X$, $\Lambda$ is of the
form $\Lambda = (L_{1,p_1},\dots,L_{n,p_n})$ as above.  Thanks to the
condition of Lemma \ref{condition2}, we see that for each $j$, replacing the
coordinate $L_{j,p_j}$ by $L_{j,p}$ for any $p$ with $p_j \leq p \leq N$ gives
a component of $X$ which is also not in $V$.  

Let $B_1$ be the set of indices $j$ for which $p_j = 1$ and let $B_{\geq 2}$
be the set of indices $j$ for which $p_j \geq 2$.  Since $V$ is not empty,
$B_{\geq 2}$ is not the empty set.  Let 
\[
F = \prod_{j \in B_{\geq 2}} \left ( \prod_{i=1}^{p_j -1} L_{j,i} \right )
\]
Clearly $F$ does not vanish on $\Lambda$.  However, any component $Q =
(L_{1,q_1} ,\dots,L_{n,q_n})$ of $X$ which is which is not in the vanishing
locus of $F$ has entries which satisfy $q_j \geq p_j$ for $j \in B_{\geq 2}$
and $q_j \geq 1 = p_j$ for $j \in B_1$.  Since $\Lambda \notin V$, we thus
have $Q \notin V$ as well. 

For the second part of the theorem, let 
\[
m_1 = X_1^{a_1}\cdots X_n^{a_n}, \ \ m_2 = X_1^{b_1}\cdots X_n^{b_n}
\]
be monomials of the same degree such that $m_1 > m_2$.  Then there is a
smallest integer $j$ such that $a_i = b_i$ for $i<j$ and $a_j > b_j$. 
Suppose $m_2 \in J$. We have to show that $m_1 \in J$.  

Without loss of generality we may assume that $b_j = a_j -1$.  Then we have
\[
a_{j+1}+ \cdots+a_n + 1 = b_{j+1}+ \cdots+b_n.
\]
Since $m_2 \in J$, Lemma \ref{not-in-J} implies that $(L_{1,b_1+1},\dots,
L_{n,b_n+1} ) \notin V$.  Hence the condition of Lemma \ref{condition3} gives
that $(L_{1,a_1+1},\dots, L_{n,a_n+1} ) \notin V$ as well, so that again by
Lemma \ref{not-in-J} we have $m_1 \in J$ as claimed.
\end{proof}

\begin{corollary}
A configuration satisfying the condition of Lemma 4.3 is arithmetically
Cohen-Macaulay.
\end{corollary}

Now we can give the answer to the question posed at the beginning of this
section, namely we identify precisely which configurations are the true
liftings of Artinian monomial ideals, and which are the true liftings of
Artinian lex-segment monomial ideals.

\begin{corollary}\label{lifting-case}
Let $V$ be a configuration of linear varieties of codimension $n$ in
$\proj{n+t-1}$.   Let $A$ be a matrix of linear forms, 
\[
A = 
\left [
\begin{array}{ccccccccc}
L_{1,1} & L_{1,2} & L_{1,3} & \dots & L_{1,N} \\
L_{2,1} & L_{2,2} & L_{2,3} & \dots & L_{2,N} \\
\vdots & \vdots & \vdots & & \vdots \\
L_{n,1} & L_{n,2} & L_{n,3} & \dots & L_{n,N}
\end{array}
\right ]
\]
such that $L_{j,i} \in K[X_j,u_1,\dots,u_t]$ and such that the entries in any
fixed row are pairwise linearly independent.  Suppose that $V$ is contained in
the complete intersection $X$ defined by the polynomials $F_j = \prod_{i=1}^N
L_{j,i}$, $1 \leq j \leq n$. Then $V$ is the lifting (in the sense of
Definition \ref{t-lifting-def}), via  $A$, of an Artinian monomial ideal in
$n$ variables if and only if the condition of Lemma \ref{condition2} holds. 
Furthermore, $V$ is the lifting, via $A$, of an Artinian lex-segment
monomial ideal if and only if in addition the condition of
Lemma \ref{condition3} holds.
\end{corollary}

\begin{proof}
The proof of Theorem \ref{lift-to-stick} (a) shows that $(F_1,\dots,F_n)$ is
a lifting of $(X_1^N,\dots,X_n^N)$, hence it is a complete intersection, and
that $X$ is reduced.  Then Theorem \ref{lift-artinian} applies.
\end{proof}

\begin{remark}\label{k-config}
The notion of a $k$-configuration of points, and related notions and
applications, have received a good deal of attention in recent years (cf.\
e.g.\ \cite{harima}, \cite{GPS}, \cite{GKS}, \cite{GHS1}, \cite{GS},
\cite{GHS2}).  Our configurations above are somewhat related to these.  In
fact, the configurations of points arising from lex-segment ideals {\em are}
special kinds of $k$-configurations, slightly more general than standard
$k$-configurations but not as general as $k$-configurations.  (It is easy to
construct a lifting matrix to produce any standard $k$-configuration, but the
need for so much collinearity prevents our obtaining an arbitrary
$k$-configuration in this way.)  

On the other hand, the configurations arising from arbitrary monomial
ideals and satisfying the condition in Lemma \ref{condition2} but not
necessarily that in Lemma \ref{condition3} are {\em not} $k$-configurations.
The key is that in the definition of a $k$-configuration, for instance
\cite{GHS1} Definitions 2.3 and 2.4, there is a strict inequality $\sigma
({\mathcal T}_i) < \alpha({\mathcal T}_{i+1})$.  The configurations
produced by Lemma \ref{condition2} only need to satisfy a weak inequality
here (although note that in addition we still require the collinearity
properties, so the weak inequality is not enough).

We will thus define a {\em weak $k$-configuration} of points in $\proj{n}$
to be a finite set of points which satisfies the definition of a
$k$-configuration as cited above, except that we require only $\sigma
({\mathcal T}_i) \leq \alpha({\mathcal T}_{i+1})$. For points in $\proj{2}$
this replacement of the strong inequality by the weak one was done in
\cite{GPS}, and the result was called a {\em weak $k$-configuration}.  
\end{remark}

In what follows we would like to extend to higher dimension the notions of
$k$-config\-urations and weak $k$-configurations of points in $\proj{n}$.  The
point of this is that our work on pseudo-lifting monomial ideals immediately
shows that such configurations exist for any allowable Hilbert function.

\begin{definition}
A {\em  $k$-configuration of dimension $d$ linear varieties} ($d \geq 0$)
is an \acm generalized stick figure of dimension $d$ whose intersection with
a general linear space of complementary dimension is a $k$-configuration.  A
{\em weak $k$-configuration of dimension $d$ linear varieties} is an \acm
generalized stick figure of dimension $d$ whose intersection with a general
linear space of complementary dimension is a weak
$k$-configuration.
\end{definition}

Note that $k$-configurations (resp.\ weak $k$-configurations) of dimension $d$
linear varieties can be produced as pseudo-liftings of Artinian lex-segment
(resp.\ non lex-segment) monomial ideals.  We can now describe the
generalized stick figures of Corollary \ref{differentiable} somewhat more
clearly.

\begin{corollary}
Let $\{c_i\}$ be a $t$-times differentiable O-sequence of dimension $t$. Then
$\{c_i\}$ is the Hilbert function of a strong $k$-configuration of dimension
$t-1$ linear varieties.  If $\{c_i\}$ is not the ``generic'' Hilbert function
then it is also the Hilbert function of a weak $k$-configuration which is not
a $k$-configuration.
\end{corollary}

\begin{proof}
The ``generic'' Hilbert function is the $t$-times differentiable Hilbert
function of dimension $t$ whose corresponding $h$-vector is 
\[
1 \ \ n \ \ \binom{n+1}{2} \ \ \binom{n+2}{3}\ \ \dots \ \ \binom{n+k}{k+1} \
\ 0.
\]
Any Artinian monomial ideal with this Hilbert function is forced to be
lex-segment.  In any other case one can find a monomial ideal which is not
lex-segment.
\end{proof}

\begin{remark}
We conclude with an observation about the behavior of liaison
under pseudo-lifting.  Again let $J$ be an Artinian monomial ideal, let $N_j$
be the maximum power of $X_j$ that occurs in a minimal generator of $J$, and
let
$N$ be the maximum of the $N_j$.  Suppose we lift $J$ to an ideal $I$ using
the lifting matrix
\[
A = 
\left [
\begin{array}{ccccccccc}
L_{1,1} & L_{1,2} & L_{1,3} & \dots & L_{1,N} \\
L_{2,1} & L_{2,2} & L_{2,3} & \dots & L_{2,N} \\
\vdots & \vdots & \vdots && \vdots \\
L_{n,1} & L_{n,2} & L_{n,3} & \dots  & L_{n,N}
\end{array}
\right ].
\]
Notice that some of the linear forms $L_{j,i}$ may not be used (if $i >
N_j$).  Form the matrix
\[
A' = \left [
\begin{array}{ccccccccc}
L_{1,N_1} & L_{1,N_1 -1} &  \dots & L_{1,1} & \dots \\
L_{2,N_2} & L_{2,N_2-1} & \dots & L_{2,1} & \dots \\
\vdots & \vdots & & \vdots \\
L_{n,N_n} & L_{n,N_n -1} & \dots  & L_{n,1} & \dots
\end{array}
\right ].
\]
Since $J$ is an Artinian monomial ideal, it includes among its generators a
power of each variable.  Hence we have contained in $J$ the complete
intersection 
\[
\tilde J = (X_1^{N_1}, \dots,X_n^{N_n}).
\]
Of course we can also lift $\tilde J$ using the matrix $A$, and we obtain a
complete intersection $\tilde I \subset I$.  The amusing fact that emerges is
that the residual ideal $[\tilde I : I]$ is the pseudo-lifting of $K$, but
using $A'$ rather than $A$!  We leave the details to the reader.  (See
\cite{migliore} for useful facts about liaison theory.)
\end{remark}



\begin{thebibliography}{999}

\bibitem{ballico-m} E.\ Ballico and J.\ Migliore, {\it Smooth Curves Whose 
Hyperplane Section is a Given Set of Points}, Comm.\ Algebra {\bf 18 (9)}
(1990), 3015--3040.\

\bibitem{macaulay} D.\ Bayer and M.\ Stillman, Macaulay: A system for
computation in algebraic geometry and commutative algebra. Source and object
code available for Unix and Macintosh computers.  Contact the authors, or
download from ftp://math.harvard.edu via anonymous ftp.

\bibitem{BM5} G.\ Bolondi and J.\ Migliore, {\em Configurations of Linear 
Projective Subvarieties}, in ``Algebraic Curves and Projective Geometry,
Proceedings (Trento, 1988),'' Lecture Notes in Mathematics, vol.\ 1389,
Springer--Verlag (1989), 19--31.

\bibitem{bruns-herzog} W.\ Bruns and J.\ Herzog, ``Cohen-Macaulay Rings,''
Cambridge studies in advanced mathematics, Cambridge University Press, 1993.

\bibitem{BE-lifting} D.\ Buchsbaum and D.\ Eisenbud, {\em Lifting modules and
a theorem on finite free resolutions}, Ring theory (Proc. Conf., Park City,
Utah, 1971), pp.\ 63--74. Academic Press,  New York, 1972. 

\bibitem{B-E} D.\ Buchsbaum and D.\ Eisenbud, {\em What makes a complex
exact?}, J.\ Algebra {\bf 25} (1973), 259--268.

\bibitem{CO} L.\ Chiantini and F.\ Orecchia, {\it Plane Sections of
Arithmetically Normal Curves in} $\proj{3}$, in ``Algebraic Curves and Projective
Geometry, Proceedings (Trento, 1988),'' Lecture Notes in Mathematics, vol.\
1389, Springer--Verlag (1989), 32--42.

\bibitem{eisenbud} D.\ Eisenbud, ``Commutative Algebra with a View toward
Algebraic Geometry,'' Graduate Texts in Mathematics 150, Springer-Verlag
(1995).

\bibitem{GGR} A.V.\ Geramita, D.\ Gregory and L.\ Roberts, {\em Monomial 
Ideals and Points in Projective Space}, J.\ Pure Appl.\ Algebra {\bf 40}
(1986), 33-62.

\bibitem{GHS1} A.V.\ Geramita, T.\ Harima and Y.S.\ Shin, {\em Extremal Point
Sets and Gorenstein Ideals}, to appear in Adv.\ Math.

\bibitem{GHS2} A.V.\ Geramita, T.\ Harima and Y.S.\ Shin, {\em An Alternative
to the Hilbert Function for the Ideal of a Finite Set of Points in
$\proj{n}$}, preprint 1998.

\bibitem{GKS} A.V.\ Geramita, H.J.\ Ko and Y.S.\ Shin, {\em The Hilbert
Function and the Minimal Free Resolution of some Gorenstein Ideals of
Codimension 4}, Comm.\ Algebra, {\bf 26} (1998), 4285-4307 

\bibitem{GMR} A.V.\ Geramita, P.\ Maroscia and L.\ Roberts, {\em The Hilbert
function of a reduced $k$-algebra}, J.\ London Math.\ Soc.\ {\bf 28} (1983),
443--452.

\bibitem{GPS} A.V.\ Geramita, M.\ Pucci and Y.S.\ Shin, {\em Smooth points of
${\mathcal G}or(T)$}, J.\ Pure Appl.\ Algebra {\bf 122} (1997), 209--241.

\bibitem{GS} A.V.\ Geramita and Y.S.\ Shin, {\em $k$-configurations in
$\proj{3}$ all have Extremal Resolutions}, to appear in J.\ Algebra. 

\bibitem{harima} T.\ Harima, {\em Some examples of unimodal Gorenstein
sequences}, J.\ Pure Appl.\ Algebra {\bf 103} (1995), 313--324.

\bibitem{hart} R.\ Hartshorne, {\em Connectedness of the Hilbert scheme},
Math.\ Inst.\ des Hautes Etudes Sci.\ {\bf 29} (1966), 261--304.

\bibitem{hart-zeuthen} R.\ Hartshorne, {\em Families of curves in ${\bf P}\sp
3$ and Zeuthen's problem}, Mem.\ Amer.\ Math.\ Soc.\ {\bf 130} No.\ 617
(1997). 

\bibitem{HTV} J.\ Herzog, N.V.\ Trung and G.\ Valla, {\em On hyperplane
sections of reduced irreducible varieties of low codimension}, J.\ Math.\
Kyoto Univ.\ {\bf 34-1} (1994), 47--72.

\bibitem{migliore} J.\ Migliore, ``Introduction to Liaison Theory and
Deficiency Modules,''  Birkh\"auser, Progress in Mathematics 165, 1998.

\bibitem{RR} L.\ Reid and L.\ Roberts, {\em Intersection points of
seminormal configurations of lines}, Algebraic $K$-theory, commutative
algebra, and algebraic geometry (Santa Margherita Ligure, 1989), 151--163,
Contemp. Math., 126, Amer.\ Math.\ Soc., Providence, RI, 1992.

\bibitem{schwartau} P.\ Schwartau, ``Liaison Addition and Monomial Ideals,''
Ph.D.\ thesis, Brandeis University (1982).

\bibitem{stanley} R.\ Stanley, {\em Hilbert functions of graded algebras},
Adv.\  Math.\ {\bf 28} (1978), 57--82.

\bibitem{taylor} D.\ Taylor, ``Ideals generated by monomials in an
$R$-sequence.'' Thesis, Chicago University (1960).

\bibitem{walter1} C.\ Walter, {\em Hyperplane sections of space curves of
small genus}, Comm.\ Algebra {\bf 22} (1994), no. 13, 5167--5174. 

\bibitem{walter2} C.\ Walter, {\em Hyperplane sections of arithmetically
Cohen-Macaulay curves}, Proc.\ Amer. Math.\ Soc.\ {\bf 123} (1995), no. 9,
2651--2656.

\end{thebibliography}
\end{document}